\documentclass[10pt, a4paper]{article}

\usepackage{fullpage}
\usepackage{amssymb}
\usepackage{amsmath}
\usepackage{amsfonts}
\usepackage{subcaption}
\usepackage{multirow}
\usepackage{tikz}
\usepackage{float}

\newtheorem{remark}{Remark}

\usepackage{xcolor}
\usepackage[colorlinks]{hyperref}
\definecolor{darkgreen}{rgb}{0.0,0.4,0.0}
\definecolor{darkred}{rgb}{0.6,0.0,0.0}
\definecolor{darkblue}{rgb}{0.0,0.0,0.5}
\definecolor{gray}{rgb}{0.5,0.5,0.5}
\definecolor{cyan}{rgb}{0.0,1.0,1.0}
\definecolor{darkcyan}{rgb}{0.0,0.5,0.5}
\definecolor{darkorange}{rgb}{0.8,0.4,0.0}
\definecolor{darkmargenta}{rgb}{0.5,0.0,0.5}
\definecolor{black}{rgb}{0.0,0.0,0.0}

\hypersetup{citecolor=blue, linkcolor=red, anchorcolor=darkblue,
filecolor=darkblue, urlcolor=darkblue}
\hypersetup{pdfauthor={Nicolae Suciu},pdftitle={Suciu}}

\begin{document}

\title{Global random walk solvers for reactive transport and biodegradation
processes in heterogeneous porous media}

\author{Nicolae Suciu$^{1, 2}\footnote{Corresponding author. \textit{Email adress: suciu@math.fau.de}}$, Florin A. Radu$^3$ \\ \\
\smallskip 
\textit{$^{1}$ Mathematics Department, Friedrich-Alexander University of Erlangen-N\"urnberg,}\\
\textit{Cauerstra{\ss}e. 11, 91058 Erlangen, Germany}\\
\textit{$^{2}$Tiberiu Popoviciu Institute of Numerical Analysis, Romanian Academy,}\\
\textit{Fantanele 57, 400320 Cluj-Napoca, Romania}\\
\textit{$^{3}$ Department of Mathematics, University of Bergen,}\\
\textit{ All\'{e}gaten 41, 5007 Bergen, Norway}}

\date{}
\maketitle

\begin{abstract}
Flow and multicomponent reactive transport in saturated/unsaturated porous media are modeled by ensembles of computational particles moving on regular lattices according to specific random walk rules. The occupation number of the lattice sites is updated with a global random walk (GRW) procedure which spreads the particles from  a lattice site with computational costs comparable to those for a single random walk step in sequential procedures. To cope with the nonlinearity and the degeneracy of the Richards equation the GRW flow solver uses linearization techniques similar to the $L$-scheme developed in finite element/volume approaches. Numerical schemes for reactive transport, coupled with the flow solver via numerical solutions for saturation and water flux, are implemented in splitting procedures. Diffusion-advection steps are solved by GRW algorithms using either biased or unbiased random walk probabilities. Since the number of particles in GRW simulations can be as large as the number of molecules involved in chemical reactions, one avoids the cumbersome problem of rescaling particle densities to approximate concentrations. Reaction steps are therefore formulated in terms of concentrations, as in deterministic approaches. The numerical convergence of the new schemes is demonstrated by comparisons with manufactured analytical solutions. Coupled flow and reactive transport problems of contaminant biodegradation described by the Monod model are further solved and the influence of flow nonlinearity/degeneracy and of the spatial heterogeneity of the medium is investigated numerically.
\end{abstract}

\begin{flushleft}

\end{flushleft}

\section{Introduction}
\label{intro}

The development of accurate numerical methods for reactive transport in natural porous media faces challenges posed by the need to solve coupled systems of nonlinear equations in conditions of highly heterogeneous properties of the medium. Besides the inherent nonlinearity of the chemical reactions, their coupling to the advection and dispersion processes influences the accuracy of the predictions based on numerical simulations. For instance, an increase of the dispersion, due to the numerical diffusion of the scheme \cite{BauseandKnabner2004,Brunneretal2012} or to inappropriate parameterizations of the dispersion coefficients \cite{Cirpkaetal1999}, leads to an artificial mixing of the species and to an overprediction of the biodegradation rate. Further, reduced accuracy of the flow solutions, in conditions of parabolic/elliptic degeneracy of the Richards equation, affects the accuracy of the transport and reaction simulations \cite{Raduetal2009}. Even more challenges may appear in the case of simulations of coupled processes in heterogeneous media, when the accuracy of the flow solutions is not only conditioned by the performance  of the numerical method but also by the parameters and correlation structure of the random-function representation of the heterogeneity \cite{Alecsaetal2020,Suciu2020}.

Most of the numerical approaches in modeling reactive transport in porous media are based on finite element \cite{Knabneretal2003,Brunneretal2012,Raduetal2008} or finite volume \cite{Cirpkaetal1999,Lie2019,Illiano2020} schemes. In these approaches, complex reaction problems with large systems of mobile and immobile species are solved by reduction schemes using linear combinations to construct non-reactive components \cite{KreutleandKnabner2005,KreutleandKnabner2007}. The numerical diffusion in advection-dominated transport problems is controlled by using higher-order schemes and adaptive time steps \cite{BauseandKnabner2004}, mixed-finite element methods using Lagrange multipliers \cite{Brunneretal2012}, or streamtube approaches \cite{Cirpkaetal1999}.

An attractive alternative to solve advection-dominated problems is provided by various particle methods, which by construction are free of numerical diffusion. Particle approaches to reactive transport were developed with
sequential particle tracking (PT) \cite{BensonandMeerschaert2008} and continuous time random walk (CTRW) \cite{Ederyetal2010} schemes, as well as with sequential biased random walk cellular automata (CA) \cite{KarapiperisandBlankleider1994,Karapiperis1995} and global procedures \cite{IzsakandLagzi2010,NagyandIzsak2011} based on global random walk (GRW) algorithms \cite{Vamosetal2003}.

Within the sequential PT and CTRW approaches a reaction takes place
when two particles representing different molecular species are
close enough to interact. The decision is made by using the
``interaction radius'', an empirical parameter related to the pore/grain
size \cite{Ederyetal2010}, or by using the probability that the two
particles will occupy the same position, constructed as a
convolution of two Gaussian probability densities
\cite{BensonandMeerschaert2008}. The latter approach was used to simulate biodegradation processes in a column experiment \cite{DingandBenson2015} and in a real field application \cite{Dingetal2017}. In these simulations, the concentrations estimated by simple binning of particle numbers were further smoothed by averaging over PT runs.

Other PT approaches estimate molecular species concentrations by smooth particle distributions provided by kernel density estimations and simulate the reaction step stochastically, with reaction probabilities depending on the distance between particles and the bandwidth of the spatial kernel \cite{Fernandez-GarciaandSanchez-Vila2011,Sole-Marietal2019}. In advection-dominated problems, the effect of the numerical diffusion induced by spatial redistribution of particles into concentrations in control volumes \cite{Cirpkaetal1999,Cuietal2014,Sanchez-VilaandFernandez-Garcia2016} may be reduced by using locally adaptive particle support volumes oriented along the local velocities, as for instance in simulations of biodegradation reactions in heterogeneous aquifers \cite{Sole-MariandFernandez-Garcia2018}. Kernel density estimations have also been used in PT modeling with reduction schemes \cite{Sole-Marietal2021} similar to those introduced in finite element approaches for multicomponent reactive systems (e.g.~\cite{KreutleandKnabner2005}).

Besides the need to map particle densities into concentrations, the coupling between transport and nonlinear reactions is another issue of concern in PT approaches \cite{Sanchez-VilaandFernandez-Garcia2016}. Accounting for dying and new created particles together with the advective-diffusive displacement of each particle can be a cumbersome task in solving complex problems for nonlinear processes. For instance, in modeling network reactions and multi-rate mass transfer in heterogeneous systems the PT method is limited to first-order reactions \cite{HenriandFernandez-Garcia2015}. In modeling biodegradation processes described by multiple-Monod expressions and geochemical reactions, an operator splitting procedure has been successfully used as an alternative to PT models based on stochastic simulations of the reaction step \cite{Cuietal2014}. In this approach, rather than adding new particles to account for mass changes, variations of reactant concentrations are accounted for by changing the mass of the existing particles. The positions of the original system of species-particles tracked at all time steps is used to map particle densities into concentrations. Further, after each PT step the reactions are modeled deterministically, by assuming that all particles in a grid cell are in contact and react to each other. This approach has been validated by comparisons with experimental data, as well as with analytical and numerical solutions, for both saturated and unsaturated flow regimes (without transitions between them). However, the method still faces some technical challenges, as for instance, representing the background oxygen concentration which fills the entire domain by adding/removing particles accounting for differences from background concentration, or the need to use a ``particle-splitting'' procedure to preserve the accuracy in cells with small numbers of particles.

In approaches based on random walk on regular lattices, there is no need to search for nearby groups of particles which can interact. A reaction takes place when a sufficient amount of particles representing reacting species meet at a lattice site \cite{IzsakandLagzi2010,Karapiperis1995}. For
instance, in the CA approach, assuming that the local chemical
equilibrium is attained on a time scale much shorter than advection
and diffusion scales, both the reaction probabilities $P_r$ and the
rate of reaction production are proportional to the physical rate
constant $k$, e.g. $P_r\sim k$ and $n_{c}\sim kn_{a}n_{b}^2$, for the reaction
$a+2b\rightarrow c$  \cite{Karapiperis1995}. The CA reaction rules use functions of occupation numbers $n_{a}$ and $n_{b}$ and reaction probabilities $P_r$ to decide whether a reaction takes place if an $a$-particle and two $b$-particles meet at the same lattice site. Even though these rules are not unique, they approximate the law of mass action, which is recovered under certain limiting conditions \cite{KarapiperisandBlankleider1994}. The species concentration is retrieved by multiplying the particle densities at lattice sites by appropriate scaling factors.

In GRW-type procedures presented in \cite{IzsakandLagzi2010,NagyandIzsak2011}, the reaction step is performed deterministically. Comparing with the PT approach with deterministic reaction step \cite{Cuietal2014} described above, we note that now the positions of all the particles are globally updated at every time step by advective and diffusive displacements. Moreover, GRW procedures use sufficiently large ensembles of particles to ensure the desired accuracy in representing concentrations by particle densities at lattice sites. As an illustration, in case of the three-component reaction considered above, the amount of product concentration in the time step $\Delta t$, $\Delta c_c=\Delta t kc_{a}c_{b}^2$, is directly converted into the amount of $c$-particles by representing the species concentrations as numbers of particles at lattice sites normalized by the number of particles used to represent a mole \cite{NagyandIzsak2011}.

The GRW approaches to reactive transport proposed in
\cite{IzsakandLagzi2010,NagyandIzsak2011} do not incorporate
advection and the diffusion step is performed by generating binomial random variables. The latter represent numbers of particles undergoing unbiased random walk jumps from a lattice site to neighboring sites. The procedure corresponds to the complete and exact GRW algorithm for diffusion introduced in \cite{Vamosetal2003} which becomes expensive in terms of computing time, with a slow linear increase for numbers of particles larger than $10^{10}$ \cite[Fig. 9]{Vamosetal2003}. Therefore, a limited number of particles was used in these diffusion-reaction simulations (e.g.~$1.5\cdot 10^5$ particles per lattice site for each molecular species in \cite{NagyandIzsak2011}) and a scaling factor was necessary to convert numbers of particles into moles. The GRW solvers for reactive transport introduced below in this article rely on the ``reduced fluctuations algorithm'' which can manage arbitrarily large numbers of particles \cite{Vamosetal2003}, also used in our recent papers \cite{Suciu2020,Suciuetal2020,Suciuetal2021}. In this way, one achieves fast and highly accurate simulations, which cannot be done with the sequential PT and biased random walk CA approaches described above. Moreover, the number of particles representing a mole can be set to the actual value of the Avogadro number, so that no scaling factors are needed to relate particle densities and species concentrations. Also, for such large ensembles of particles the GRW simulations are self-averaging \cite{Vamosetal2003}, i.e. there is no need to perform averages over runs \cite{DingandBenson2015,Dingetal2017} or kernel density estimations \cite{Fernandez-GarciaandSanchez-Vila2011,Sole-Marietal2019} to obtain smooth concentrations as in PT procedures.

In the present study we propose two-dimensional GRW solvers for nonlinear reactive transport, with a special emphasis on aerobic biodegradation. The solvers consist of explicit numerical schemes which solve the transport and the reaction steps by an operator splitting procedure. The flow velocity and the water content are provided by the GRW solver presented in \cite{Suciuetal2021}, where the treatment of the nonlinear terms in the Richards equation is based on the $L$-scheme linearization approach \cite{Popetal2004, ListandRadu2016}. We consider here only the one-way coupling of Richards and transport equations in which the concentrations do not influence the flow.

With the unbiased GRW, the transport is modeled by updating the occupation numbers of the computational lattice with a global procedure which relocates all the particles from a lattice site in a single numerical procedure. First, all the particles undergo advective displacements, according to the local velocity evaluated at the lattice site, then they are distributed to neighboring sites by computing the number of unbiased random walk jumps in the spatial directions of the lattice. Considering the movement of a single particle, we note that it is similar to the PT algorithm, with the difference that the advective displacement is discrete and the diffusive one is modeled as unbiased random walk, rather than by realizations of the Wiener process (normally distributed random numbers) as in PT approaches. In fact, the GRW algorithm is a superposition of PT schemes (or weak solutions of It\^{o} equation) on a regular lattice \cite[Sect. 3.3]{Suciu2019}.

Instead, if the jump probabilities are biased, the bias being proportional to the local velocity, the random walk models both the advective and diffusive displacements. This results in a biased global random walk (BGRW) algorithm which generalizes the CA approach of \cite{KarapiperisandBlankleider1994} as a global procedure. In the BGRW algorithm the particles are allowed to jump only to first-neighbor sites. In this way, one avoids jumps over several lattice sites with significant differences in local velocity and diffusion coefficients which produce overshooting errors. Thus, the BGRW algorithm is more accurate than the unbiased GRW, which instead, by using a coarser spatial discretization is faster.

In both unbiased GRW and BGRW solvers, the reaction step is solved deterministically by evaluating the reaction terms in the coupled system of reactive transport equations. In case of saturated flows, the nonlinear reaction terms are evaluated with an explicit linearization by computing them with concentration values from the previous time steps of the numerical scheme. Since numerical tests indicate that the explicit linearization is no longer accurate for variable water content, we propose an iterative approach for reactive transport in variably saturated flows similar to the GRW $L$-schemes introduced in \cite{Suciuetal2020,Suciuetal2021}.

The paper is organized as follows. In Section~\ref{sec:model} we introduce the governing equations for coupled flow and reactive transport with Monod reactions. Section~\ref{sec:solver} contains the presentation of the GRW solvers for multicomponent reactive transport. The solvers are further evaluated in Section~\ref{sec:verif} by comparisons with analytical solutions and numerical estimations of orders of convergence for saturated and unsaturated flow regimes, as well as in case of degenerate Richards equation describing the transition between unsaturated and saturated flows. In Section~\ref{sec:bioSat} the new solvers are tested on benchmark problems for biodegradation processes governed by Monod type models in saturated flows often used in studies published in the past. Further, the reactive transport solver will be coupled to Richards equation with randomly distributed parameters and the influence of the variable saturation and of the spatial heterogeneity of the soil/aquifer system will be investigated through numerical experiments in Section~\ref{sec:deg_Rich}. Section~\ref{sec:concl} contains the conclusions of this study. The GRW/BGRW codes implemented in Matlab used in this study are stored in the Git repository \href{https://github.com/PMFlow/MonodReactions}{MonodReactions} \cite{SuciuandRadu2021}.

\section{A simple model for reactive transport and biodegradation}
\label{sec:model}

The microbial biodegradation is a natural process which facilitates the in situ bioremediation of contaminated soils and groundwater systems, either as an intrinsic process or as part of an engineered remediation technology  \cite{Wiedemeieretal1999}. The process occurs when both an electron acceptor (e.g. oxygen) and a microbian population are available. The contaminant to be degraded (e.g. benzene) acts as an electron donor and as a carbon source utilized by microorganisms \cite{BauseandKnabner2004}. Since the availability of the three elements depends on the water flow and the diffusion/dispersion processes, reliable models should consider the coupling of the reaction system to the flow and transport equations \cite{Radu2004}.

We consider a two-dimensional setting with an open domain $\Omega$, with boundary $\Gamma$, and a time interval $(0,T]$. Taking into account both the saturated and the vadose zone, the water flux is given by Darcy's law
\begin{equation}\label{eq:Q}
\mathbf{q}=-K(\theta(\psi)\nabla(\psi+z),
\end{equation}
where $z$ is the height against the gravitational direction. The pressure head $\psi(x,z,t)$ solves the Richards equation
\begin{equation}\label{eq:Richards}
\frac{\partial}{\partial t}\theta(\psi)-\nabla\cdot\left[K(\theta(\psi)\nabla(\psi+z)\right]=0,
\end{equation}
where $\theta$ is the water content and $K$ stands for the hydraulic conductivity of the medium. Equation (\ref{eq:Richards}) is closed by material laws expressing dependencies between $K$, $\theta$, and $\psi$, based on experimental results. Further, we consider three molecular species, of which the first two are mobile and the third is immobile. The evolution of the corresponding concentrations  $c_\nu$, $\nu=1,2,3$ is governed by the following system of coupled nonlinear equations
\begin{align}
&\frac{\partial}{\partial t}\left(\theta c_\nu\right)-\nabla\cdot\left(D_\nu\nabla c_\nu -\mathbf{q}c_\nu\right)=R_\nu, \;\; \nu=1,2,\label{eq:Mobile}\\
&\frac{\partial}{\partial t}c_3=R_3\left(1-\gamma\frac{c_3}{c_{3_{max}}}\right)-k_{d}c_{3}.\label{eq:Immobile}
\end{align}
The constant $k_d$ in Eq.~(\ref{eq:Immobile}) is the decay rate of the immobile component. The factor of the reaction rate $R_3$ may be used to model the growth limitation of $c_3$ (due for instance to limited pore space available), $c_{3_{max}}$ is the maximum immobile concentration, and the coefficient $\gamma$ is set to one in the model with growth limitation and to zero otherwise \cite{Radu2004,BauseandKnabner2004}.  In a basic scenario for aerobic biodegradation \cite{BauseandKnabner2004,Brunneretal2012,Cirpkaetal1999,Radu2004}, $c_1$ represents the concentration of the electron donor, $c_2$ corresponds to the electron acceptor (oxygen), $c_3$ is the biomass concentration, and the growth of the biomass is assumed to follow a double Monod kinetics. The reaction term of Eq.~(\ref{eq:Immobile}) is expresses as $R_3=Y\mu$, $Y$ is the microbial yield coefficient per unit electron donor consumed, and the Monod term is defined by
\begin{equation}\label{eq:Monod}
\mu=\mu_{max}\frac{c_1}{M_1+c_1}\frac{c_2}{M_2+c_2}c_3,
\end{equation}
where $\mu_{max}$ is the maximum growth rate and $M_1$, $M_2$ are Monod constants. The reaction terms for the mobile species in Eqs.~(\ref{eq:Mobile}) are expressed as $R_1=-\theta\alpha_1\mu$ and $R_2=-\theta\alpha_2\mu$, where $\alpha_1$, $\alpha_2$ are stoichiometric constants.

As compared with more complex biodegradation models \cite{BauseandKnabner2004,Cuietal2014,Radu2004}, the model presented above is simplified in that only the degradation of a single contaminant is considered. Also, the inhibition terms which yield a slower microbian growth and the sorption terms are neglected. Instead, the coupling of the basic model process of aerobic degradation described by (\ref{eq:Mobile})-(\ref{eq:Monod}) with the flow model (\ref{eq:Q})-(\ref{eq:Richards}) allows investigations on the role of the variable saturation and of the heterogeneity of the soil-aquifer system.

\section{Splitting algorithms for multicomponent reactive transport}
\label{sec:solver}
In the biodegradation model presented above, the advection-dispersion-reaction problem is coupled in one way with the flow problem through the solution $(\theta, \mathbf{q})$ of the flow equations (\ref{eq:Q})-(\ref{eq:Richards}). The system of equations (\ref{eq:Q})-(\ref{eq:Monod}) can be solved with methods based on linearization procedures developed for fully coupled flow and transport equations \cite{Illiano2020,Knabneretal2003,Suciuetal2021} after replacing the general constitutive law for the nonlinear water content $\theta(\psi,c)$ by the simpler law $\theta(\psi)$, corresponding to the one-way coupling. The resulting scheme has to be completed by the computation of the coupled nonlinear reaction terms at every time step and iteration in the transport solver. With these modifications in the scheme for fully coupled flow and transport from \cite{Suciuetal2021} one obtains an iterative GRW solver for multicomponent reactive transport. We note that the operator splitting PT method for reactive transport in variably saturated media proposed in \cite{Cuietal2014} also uses an iterative Newton method to solve the nonlinear system of ordinary differential equations of chemical kinetics. In the saturated flow regime, $\theta$ no longer depends on the pressure $\psi$ and, in the examples presented in this article, it will be taken as a constant parameter. For given flow velocity, the GRW solution can be constructed, similarly to the solution for reaction-diffusion problems from \cite{NagyandIzsak2011}, as an explicit non-iterative scheme which evaluates nonlinear reaction terms from concentrations of mobile species obtained after the advection-dispersion step and immobile species concentrations from the previous time step (i.e., by an explicit linearization \cite{Haverkampetal1977}).

\subsection{Non-iterative schemes for reactive transport in saturated media}
\label{sec:noniter}
\subsubsection{Biased GRW algorithm for the transport step}
\label{sec:bgrw}

For constant $\theta$, the transport equation (\ref{eq:Mobile}) for one mobile concentration $c$, with a source/sink term $f$ independent of $c$ added in the right-hand side, can be equivalently written as
\begin{equation}\label{eq:EqivMobile}
\frac{\partial c}{\partial t}-\frac{1}{\theta}\nabla\cdot\left(D\nabla c -\mathbf{q}c\right)=\frac{R}{\theta}+\frac{f}{\theta}.
\end{equation}
The new term $f$ of the transport equation will be used in code verification tests presented in Section~\ref{sec:verif} to describe the result of applying the differential operators of the transport equation to a prescribed analytical solution.

The flow velocity $\mathbf{q}$ is computed with the GRW algorithm for flows in saturated/unsaturated porous media presented in \cite[Sect. 4.1]{Suciuetal2021}, with velocity on the boundaries extended from first-neighbor interior lattice sites or approximated by forward finite differences, for the particular case of constant $\theta$. The transport step, described by Eq.~(\ref{eq:EqivMobile}) without the reaction term, can be computed with explicit GRW algorithms. In the following, a two-dimensional BGRW which simulates the advection through biased jump probabilities \cite{KarapiperisandBlankleider1994,Suciu2019} is derived from the forward-time centered-space finite difference approximation of the transport equation. We assume a diagonal diffusion tensor with constant components $D_1$ and $D_2$ and by $U$ and $V$ we denote the components of the Darcy velocity $\mathbf{q}$ along the horizontal axis $x$ and the vertical axis $z$, by $\Delta t$ the time step, and by $\Delta x$ and $\Delta z$ the spatial steps. With these, the explicit scheme reads
\begin{align}\label{eq:FD}
c_{i,j,k+1}=c_{i,j,k}&-\frac{\Delta t}{2\theta\Delta x}\left(U_{i+1,j,k}c_{i+1,j,k}-U_{i-1,j,k}c_{i-1,j,k}\right)-\frac{\Delta t}{2\theta\Delta z}\left(V_{i,j+1,k}c_{i,j+1,k}-V_{i,j-1,k}c_{i,j-1,k}\right)\nonumber\\
&+\frac{D_1\Delta t}{\theta\Delta x ^2}\left(c_{i+1,j,k}-2c_{i,j,k}+c_{i-1,j,k}\right)+\frac{D_2\Delta t}{\theta\Delta z ^2}\left(c_{i,j+1,k}-2c_{i,j,k}+c_{i,j-1,k}\right)+\frac{f_{i,j,k}\Delta t}{\theta}\nonumber\\
=c_{i,j,k}&-\left(\frac{2D_1\Delta t}{\theta\Delta x ^2}+\frac{2D_2\Delta t}{\theta\Delta z ^2}\right)c_{i,j,k}\nonumber\\
&+\left(\frac{D_1\Delta t}{\theta\Delta x ^2}-\frac{\Delta t}{2\theta\Delta x}U_{i+1,j,k}\right)c_{i+1,j,k}+\left(\frac{D_1\Delta t}{\theta\Delta x ^2}+\frac{\Delta t}{2\theta\Delta x}U_{i-1,j,k}\right)c_{i-1,j,k}\nonumber\\
&+\left(\frac{D_2\Delta t}{\theta\Delta z ^2}-\frac{\Delta t}{2\theta\Delta z}V_{i,j+1,k}\right)c_{i,j+1,k}+\left(\frac{D_2\Delta t}{\theta\Delta z ^2}+\frac{\Delta t}{2\theta\Delta z}V_{i,j-1,k}\right)c_{i,j-1,k}+\frac{f_{i,j,k}\Delta t}{\theta}.
\end{align}

As usual in GRW algorithms, we simulate the random walk of $\mathcal{N}$  computational particles to approximate the concentration by the density of the number of particles at lattice sites $(i\Delta x,j\Delta z)$ and time points $k\Delta t$, $c_{i,j,k}\approx n_{i,j,k}/\mathcal{N}$. By using the finite difference scheme (\ref{eq:FD}) and the dimensionless parameters
\begin{equation}\label{eq:ParBGRW1}
r_{x}=\frac{2D_1\Delta t}{\theta\Delta x ^2},\;\; r_{z}=\frac{2D_2\Delta t}{\theta\Delta z ^2},\;\; u_{i\pm 1,j,k}=\frac{\Delta t}{\theta\Delta x}U_{i\pm 1,j,k},\;\; v_{i,j\pm 1,k}=\frac{\Delta t}{\theta\Delta z}V_{i,j\pm 1,k},
\end{equation}
we obtain
\begin{align}\label{eq:BGRW1}
n_{i,j,k+1}=
&\left[1-\left(r_{x}+r_{z}\right)\right]n_{i,j,k}\nonumber\\
& + \frac{1}{2}\left(r_{x}-u_{i+ 1,j,k}\right)n_{i+1,j,k}+ \frac{1}{2}\left(r_{x}+u_{i-1,j,k}\right)n_{i-1,j,k}\nonumber\\
& + \frac{1}{2}\left(r_{z}-v_{i,j+ 1,k}\right)n_{i,j+1,k}+ \frac{1}{2}\left(r_{z}+v_{i,j-1,k}\right)n_{i,j-1,k} +\left\lfloor \mathcal{N}\frac{f_{i,j,k}\Delta t}{\theta}\right\rfloor,
\end{align}
where $\lfloor\cdot\rfloor$ denotes the floor function.

The contributions to $n_{i,j,k+1}$ in Eq.~(\ref{eq:BGRW1}) are obtained with the BGRW algorithm
\begin{equation}\label{eq:BGRW2}
n_{l,m,k}=\delta n_{l,m\;\mid\; l,m,k} + \delta n_{l-1,m\;\mid\; l,m,k}+\delta n_{l+1,m\;\mid\; l,m,k} + \delta n_{l,m-1\;\mid\; l,m,k}+\delta n_{l,m+1\;\mid\; l,m,k},
\end{equation}
where, for consistency with the finite difference scheme~(\ref{eq:BGRW1}), the quantities $\delta n$ verify in the mean
\begin{align}\label{eq:BGRW3}
&\overline{\delta n_{l,m\;\mid\; l,m,k}}=\left[1-\left(r_{x}+r_{z}\right)\right]\overline{n_{,lm,k}},\nonumber\\
&\overline{\delta n_{l\pm 1,m\;\mid\; l,m,k}}=\frac{1}{2}(r_{x}\pm u_{l,m,k})\overline{n_{l,m,k}},\nonumber\\
&\overline{\delta n_{l,m\pm 1\;\mid\; l,m,k}}=\frac{1}{2}(r_{z}\pm
v_{l,m,k})\overline{n_{l,m,k}}.
\end{align}
The binomial random variables $\delta n$ used in the BGRW algorithm are approximated by using the unaveraged relations (\ref{eq:BGRW3}), summing up the remainders of multiplication by $r$ and of the floor function $\left\lfloor \mathcal{N}f/\theta\right\rfloor$, and allocating one particle to the lattice site where the sum reaches the unity (see also \cite{Suciu2019,Suciuetal2020}). By giving up the particle's indivisibility and using the un-averaged relations (\ref{eq:BGRW3}) one obtains a deterministic BGRW algorithm. However, throughout the present study we shall use only the randomized implementation of the BGRW algorithm.

\begin{remark}
\label{rem:part_nr}
(B)GRW simulations with the ``reduced fluctuations algorithm'' model the collective behavior of arbitrarily large numbers of particles at computational costs comparable to those for a single random walk or particle tracking step in sequential procedures. This feature of the GRW algorithms leads to smooth solutions and to easier implementations of the reaction schemes, by using as many particles as the number of molecules involved in reactions.
\end{remark}

As follows from (\ref{eq:BGRW3}), the BGRW algorithm is subject to the
following restrictions
\begin{equation}
r_{x}+r_{z}\leq1,\;\left\vert u_{l,m,k}\right\vert \leq
r_{x},\;\left\vert v_{l,m,k}\right\vert \leq r_{z}.\label{eq:RestrBGRW}
\end{equation}

From the definition (\ref{eq:ParBGRW1}) of the BGRW parameters and the constraints (\ref{eq:RestrBGRW}) it follows that the allowable value of the local P\'{e}clet number is limited by $\text{P\'{e}}\le 2$ (see \cite[Remark 4]{Suciuetal2020}).

\subsubsection{Unbiased GRW algorithm for the transport step}
\label{sec:ugrw}

The unbiased two-dimensional GRW algorithm for the transport step in Eq.~(\ref{eq:EqivMobile}) is obtained by globally moving groups of particles according to the rule
\begin{align}\label{eq:GRW1}
n_{l,m,k}=&\;\delta n_{l+u_{l,m,k},m+v_{l,m,k}\;\mid\; l,m,k}\\
&+ \delta n_{l+u_{l,m,k}+d,m+v_{l,m,k}\;\mid\; l,m,k}+\delta n_{l+u_{l,m,k}-d,m+v_{l,m,k}\;\mid\; l,m,k}\nonumber \\
&+ \delta n_{l+u_{l,m,k},m+v_{l,m,k}+d\;\mid\; l,m,k}+\delta n_{l+u_{l,m,k},m+v_{l,m,k}-d\;\mid\; l,m,k},\nonumber
\end{align}
where $d$ is a constant amplitude of diffusion jumps and the dimensionless variables $r_x$, $r_z$, $u$ and $v$ are defined similarly to (\ref{eq:ParBGRW1}) by
\begin{equation}\label{eq:GRW2}
r_x=\frac{2D_1\Delta t}{\theta(d\Delta x )^2},\;\; r_z=\frac{2D_2\Delta t}{\theta(d\Delta z )^2},\;\; u_{l,m,k}=\left\lfloor\frac{\Delta t}{\theta\Delta x}U_{l,m,k}+0.5\right\rfloor,\;\; v_{l,m,k}=\left\lfloor\frac{\Delta t}{\theta\Delta z}V_{l,m,k}+0.5\right\rfloor.
\end{equation}

The particle distribution is updated at every time step by
\begin{equation}
n_{i,j,k+1}=\delta n_{i,j,k}+\sum_{l\neq i,m\neq j}\delta n_{i,j\;\mid\; l,m,k} +\left\lfloor \mathcal{N}\frac{f_{i,j,k}\Delta t}{\theta}\right\rfloor.\label{eq:GRW3}
\end{equation}

The averages over GRW runs of the terms from (\ref{eq:GRW1}) are now
related by
\begin{align}
&\overline{\delta n_{l+u_{l,m,k},m+v_{l,m,k}\;\mid\; l,m,k}}=\left[1-\left(r_{x}+r_{z}\right)\right]\overline{n_{l,m,k}},\nonumber\\
&\overline{\delta n_{l+u_{l,m,k}\pm d,m+v_{l,m,k}\;\mid\; l,m,k}}=\frac{r_{x}}{2}\hspace{0.1cm}\overline{n_{l,m,k}},\nonumber\\
&\overline{\delta n_{l+u_{l,m,k},m+v_{l,m,k}\pm d\;\mid\; l,m,k}}=\frac{r_{z}}{2}\hspace{0.1cm}\overline{n_{l,m,k}}.\label{eq:GRW4}
\end{align}
Comparing with the BGRW relations (\ref{eq:BGRW3}), we remark that (\ref{eq:GRW4}) defines unbiased jump probabilities $r_x/2$ and $r_y/2$ in the two spatial directions.

The binomial random variables $\delta n$ used in the unbiased GRW algorithm are approximated by the procedure used for the BGRW algorithm presented in the previous subsection. For fixed space steps, the time step size is chosen such that the dimensionless parameters $u_{l,m,k}$ and $v_{l,m,k}$ take integer values larger than unity which ensure the desired resolution of the velocity components \cite[Sect. 3.3.2.1]{Suciu2019}. Further, the amplitude of the jumps $d$ is chosen such that the jump probabilities verify the constraint $r_{x}+r_{z}\leq 1$, imposed by the first relation (\ref{eq:GRW4}).

\subsubsection{Deterministic reaction step}
\label{sect:react}

Proceeding as in \cite{NagyandIzsak2011}, we express the species concentration in moles per lattice site. The classical unit of moles per liter can simply be obtained by dividing the number of molecules by a fixed volume in liters and by the number of molecules in a mole, as done, for instance, in the simplified molecular dynamics approach to chemical kinetics proposed in \cite{HighamandKloeden2021}.  Since the (B)GRW implementations are not sensible to the increase of the number of particles (Remark~\ref{rem:part_nr}) we represent a mole of substance by the Avogadro number $N_A\approx 6.022\cdot 10^{23}$, hereafter rounded to $\mathcal{N}=10^{24}$. The concentrations of the mobile molecular species $\nu=1, 2$ are computed as $c_{\nu,i,j,k+1/2}=n_{\nu,i,j,k+1}/\mathcal{N}$, where $n_{\nu,i,j,k+1}$ is the distribution of the computational particles after the execution of the transport step obtained with either the BGRW algorithm (\ref{eq:BGRW1}) or the unbiased GRW algorithm (\ref{eq:GRW3}). With this representation of the species concentrations, the reaction step in (\ref{eq:EqivMobile}), for the particular case of the Monod model (\ref{eq:Mobile})-(\ref{eq:Monod}), will be performed deterministically (see also \cite{NagyandIzsak2011}) according to
\begin{equation}\label{eq:React}
c_{\nu,i,j,k+1}=c_{\nu,i,j,k+1/2}+\frac{\Delta t}{\theta}R(c_{1,i,j,k+1/2},c_{2,i,j,k+1/2},c_{3,i,j,k}),
\end{equation}
where $c_{3,i,j,k}$ is the concentration of the immobile species at the previous time step.

Following the reaction step, the numbers of particles representing mobile species are updated by $n_{\nu,i,j,k+1}=c_{\nu,i,j,k+1}\mathcal{N}$, $\nu=1,2$ and the splitting procedure is resumed by executing the new flow and transport steps.

\subsection{Iterative scheme for reactive transport in variably saturated media}
\label{sec:iter}

Since in case of variably saturated soils the unbiased GRW transport solver often requires extremely fine discretizations to obtain acceptable resolution of the discrete velocity components given by (\ref{eq:GRW2}) \cite[Sect. 5.2.3]{Suciuetal2021} we shall consider in the following only the biased algorithm. The iterative BGRW $L$-scheme is obtained from a backward-time centered-space finite difference scheme for the transport equation for mobile species, Eq.~(\ref{eq:Mobile}) without reaction term and with a source term $f$ added in the right-hand side, by adding the stabilization term $L(c_{i,j,k}^{s+1}-c_{i,j,k}^{s})$, where $L$ is a dimensionless constant parameter and $s$, $s=1,2,\ldots $, is the iteration index. Following the same procedure as in Section~\ref{sec:bgrw} one obtains
\begin{align}\label{eq:FT3}
n_{i,j,k}^{s+1}
=&\left[1-\left(r_{x}+r_{z}\right)\right]n_{i,j,k}^{s}\nonumber\\
& + \frac{1}{2}\left(r_{x}-u_{i+ 1,j,k}^{s}\right)n_{i+1,j,k}^{s}+ \frac{1}{2}\left(r_{x}+u_{i-1,j,k}^{s}\right)n_{i-1,j,k}^{s}\nonumber\\
& + \frac{1}{2}\left(r_{z}-v_{i,j+ 1,k}^{s}\right)n_{i,j+1,k}^{s}+ \frac{1}{2}\left(r_{z}+v_{i,j-1,k}^{s}\right)n_{i,j-1,k}^{s} +\left\lfloor g^{s}_{i,j,k}\right\rfloor,
\end{align}
where $g^{s}_{i,j,k}=\mathcal{N}f_{i,j,k}\Delta t/L-\left[\theta(\psi_{i,j,k}^{s})n_{i,j,k}^{s}
-\theta(\psi_{i,j,k-1})n_{i,j,k-1}\right]/L$, with pressure head $\psi$ provided by the flow solver (see also \cite[Sect. 4.2.1]{Suciuetal2021}).

The parameters of the BGRW $L$-scheme are defined by
\begin{equation}\label{eq:ParBGRW}
r_{x}=\frac{2D_1\Delta t}{L\Delta x ^2},\;\; r_{z}=\frac{2D_2\Delta t}{L\Delta z ^2},\;\; u_{i\pm 1,j,k}^{s}=\frac{\Delta t}{L\Delta x}U_{i\pm 1,j,k}^{s},\;\; v_{i,j\pm 1,k}^{s}=\frac{\Delta t}{L\Delta z}V_{i,j\pm 1,k}^{s}.
\end{equation}
The BGRW $L$-scheme is completed by the same relations (\ref{eq:BGRW2})-(\ref{eq:RestrBGRW}) as in case of non-iterative BGRW scheme from Section~\ref{sec:bgrw}, which now have to be verified at every iteration $s$.
The iterations start with $c_{i,j,k}^{1}=c_{i,j,k-1}$ and are stopped when the discrete $L^2$ norms of the numerical solution $c_{k}^{s}=(c_{i,j,k}^{s})$ verify
\begin{equation}\label{eq:5}
\|c_{k}^{s}-c_{k}^{s-1}\|\leq \varepsilon_a + \varepsilon_r\|c_{k}^{s}\|
\end{equation}
for some given tolerances $\varepsilon_a$ and $\varepsilon_r$.

This BGRW $L$-scheme is embedded into an operator splitting procedure alternating nonreactive transport for all mobile species $c_\nu$ and reaction steps according to Eq.~(\ref{eq:React}) at every iteration level $s$.

\section{Code verification tests}
\label{sec:verif}

In the absence of exact solutions to boundary value problems for advective-dispersive transport coupled with Monod models, the codes can be verified by comparisons with analytical solutions to less complex problems, see e.g. \cite{Cuietal2014}, where one uses solutions for linear degradation and saturated flows with constant velocity. Using manufactured analytical solutions, instead, allows code verification tests for nontrivial situations such as flows governed by degenerate Richards equation and one-way coupling with nonlinear reactions \cite{Raduetal2009} or fully coupled one-dimensional flow and transport with transition from unsaturated to saturated conditions \cite{Suciuetal2020}. Here, we adopt this approach to test the new schemes for problems of increased complexity: constant flow and nonlinear reactions, saturated flow coupled with Monod biodegradation reactions, and space-time variable flow with transition from unsaturated to saturated regime coupled with Monod reactions.

\subsection{Constant velocity and bimolecular nonlinear reactions}
\label{sec:const_flow}

To verify the correctness of the transport and reaction steps in the splitting procedure we consider the system of coupled equations (\ref{eq:Mobile}) with given parameters $U=0$, $V=-1$, $D_1=D_2=D=0.1$, $\theta=1$, and reaction terms given by $R_1=-c_{1}c_{2}^{2}$ and $R_2=-2c_{1}c_{2}^{2}$. The problem is solved in the domain $\Omega=(0,2)\times(0,3)$ and the final simulation time is $T=1$. The manufactured analytical solutions
\begin{align}
&c_{1m}(x,z,t) = x(2-x)z^3\exp(-0.1t)/27,\label{eq:c1}\\
&c_{2m}(x,z,t) = (x-1)^{2}z^2\exp(-0.1t)/9,\label{eq:c2}
\end{align}
induce supplementary source terms in the right hand side of Eqs.~(\ref{eq:Mobile}) and define the initial and boundary conditions \cite{Raduetal2008,Brunneretal2012}.

Numerical solutions obtained with the non-iterative BGRW algorithm presented in Section~\ref{sec:bgrw} are computed on regular lattices successively refined by halving five times the initial step size $\Delta x = \Delta z = 0.2$. The accuracy of the numerical solutions at the final time $t=T$ is quantified by discrete $L^2$ norms and the estimated order of convergence (EOC) is given by the slope of decreasing error norms in logarithmic scale \cite{Alecsaetal2020,Suciuetal2021}.

\begin{table}[!ht]
\caption{Convergence of the non-iterative BGRW scheme for\\ constant velocity and nonlinear reactive transport.}
\label{table_BGRWconv}
\begin{tabular}{ c c c c c c c  c c c}
\hline
   & $\Delta x$ & $\|c_1-c_{1m}\|$ & EOC  & $\|c_2-c_{2m}\|$ & EOC \\
  \hline
   &  2.00e-1  & 3.53e-03  & -- & 4.91e-03 & -- &\\
   &  1.00e-1  & 8.64e-04  & 2.03 & 1.12e-03 & 2.13 &\\
   &  5.00e-2  & 2.12e-04  & 2.02 & 2.67e-04 & 2.07 &\\
   &  2.50e-2  & 5.30e-05  & 2.00 & 6.50e-05 & 2.04 &\\
   &  1.25e-2  & 1.32e-05  & 2.00 & 1.60e-05 & 2.02 &\\
   &  6.25e-2  & 3.30e-06  & 2.00 & 3.99e-06 & 2.01 &\\
  \hline
\end{tabular}
\end{table}
The results of the verification test presented in Table~\ref{table_BGRWconv} show a very good accuracy of the BGRW solver for reactive transport and indicate the expected convergence of order 2. A similar test was done for the one-dimensional version of the BGRW solver, with manufactured solution similar to those for the two-dimensional case, where the factors containing the $x$-variable were removed. In this case, it is found that the numerical solutions converge with the order 2 as well, but the $L^2$ errors are generally one order of magnitude larger (similarly to the one-dimensional solutions for the coupled flow and transport problem presented in \cite[Tables 10 and 11]{Suciuetal2020}).

For the coarsest lattice used in this test ($\Delta x=0.2$) the local P\'{e}clet number $\text{P\'{e}}=|V|\Delta x/D=2$ takes the maximum allowable value (see Section~\ref{sec:bgrw}). Hence, as $\Delta x$ decreases $\text{P\'{e}}<2$, the relations (\ref{eq:RestrBGRW}) are verified, and the BGRW algorithm is free of numerical diffusion \cite[Table 9]{Suciuetal2021}.
For smaller $D$ the advection dominates and $\Delta x$ has to be decreased to fulfil the condition $\text{P\'{e}}\le 2$, which would slow down the BGRW computations. The alternative is to use the unbiased GRW solver presented Section~\ref{sec:ugrw} which is not constrained by a maximum $\text{P\'{e}}$. By varying the dimensionless velocity $v=\left\lfloor\Delta t |V|/\Delta z+0.5\right\rfloor$ and the amplitude of diffusion jumps $d$ we have two degrees of freedom which can be used to accommodate different values of the velocity and diffusion coefficient. On the other side, while boundary conditions can be easily formulated for $v=0$ and $d=1$ \cite{Vamosetal2003}, if $v>1$ and $d>1$ there are groups of particles jumping across the boundary from interior sites (e.g.,~$\delta n_{i,j-v-d\;\mid\; i,j,k}$ for $j-v-d<1$ in the test problem considered here), thus inevitably inducing numerical errors. In such cases, very fine lattices would be necessary to reduce the approximation errors. Note that this issue is avoided in simulations of transport in unbounded domains \cite[Sect. 6.3]{Suciuetal2020}, which renders the unbiased GRW very efficient in large scale simulations of transport in groundwater \cite{Suciu2019}.

The smallest errors of the GRW solutions are obviously obtained if $v=1$ and $d=1$. For the setup of the present test problem this choice is not possible because if $v=1$ the amplitude of the diffusion jumps in (\ref{eq:GRW2}) should be at least  $d=2$ to comply with the restriction $r_x+r_z\le 1$, for $\Delta x=0.2$, and even larger $d$ is needed for smaller $\Delta x$. Instead, for strongly advection-dominated problems with $D=10^{-4}$, as those considered in the next section, the choice $v=1,\; d=1$ is possible for all $\Delta x$ considered in Table~\ref{table_BGRWconv}. The results of the biased GRW solutions, with the new dispersion coefficient and letting the other parameters unchanged, are presented in Table~\ref{table_GRWconv}.
\begin{table}[!ht]
\caption{Convergence of the non-iterative unbiased GRW scheme\\ for constant velocity and nonlinear reactive transport.}
\label{table_GRWconv}
\begin{tabular}{ c c c c c c c  c c c}
\hline
   & $\Delta x$ & $\|c_1-c_{1m}\|$ & EOC  & $\|c_2-c_{2m}\|$ & EOC \\
  \hline
   &  2.00e-1  & 1.69e-01  & --   & 7.72e-02 & --   &\\
   &  1.00e-1  & 8.43e-02  & 1.00 & 3.01e-02 & 1.36 &\\
   &  5.00e-2  & 4.13e-02  & 1.03 & 1.27e-02 & 1.24 &\\
   &  2.50e-2  & 2.07e-02  & 0.99 & 5.96e-03 & 1.10 &\\
   &  1.25e-2  & 1.04e-02  & 1.00 & 2.91e-03 & 1.04 &\\
   &  6.25e-2  & 5.14e-03  & 1.01 & 1.42e-03 & 1.03 &\\
  \hline
\end{tabular}
\end{table}
The larger error norms and the slower convergence of order 1 shown in the table can be attributed to the inherent errors of the GRW solver due to jumps over the lower boundary of the computational domain.

\subsection{Saturated flow and Monod reactions}
\label{sec:sat_flow}

We consider in the following the same computational domain and final time as in Section~\ref{sec:const_flow} above. A more challenging test problem is formulated by considering an unsteady saturated flow with pressure head given analytically by
\begin{equation}
\psi_m(x,z,t) = tx(2-x)z(3-z),\label{eq:p_sat}
\end{equation}
together with a constant water content $\theta=0.3$, constant hydraulic conductivity $K=0.05$, and constant dispersion coefficients $D_1=D_2=0.025$. The nonlinear reaction is governed by the double Monod model (\ref{eq:Monod}) with parameters $\mu_{max}=10^{-3}$, $M_1=2$, $M_2=0.2$, stoichiometric coefficients $\alpha_1=1$,  $\alpha_2=3$, analytical solutions for the election donor and election acceptor concentrations verifying Eqs.~(\ref{eq:Mobile})-(\ref{eq:Immobile}) given by (\ref{eq:c1})-(\ref{eq:c2}), and constant biomass concentration $c_3=1$.

The coupled system of equations (\ref{eq:Richards})-(\ref{eq:Mobile}) with source terms in the right-hand side, computed analytically by differentiating the solutions (\ref{eq:c1})-(\ref{eq:p_sat}), is solved with the flow GRW $L$-scheme from \cite[Sect. 4.1]{Suciuetal2021} coupled with the non-iterative BGRW scheme introduced in Section~\ref{sec:bgrw}. The flow solutions are approximated, with parameter $L=1$ and tolerances $\varepsilon_a=\varepsilon_r=5\cdot 10^{-7}$, by performing a large number of iterations of the flow scheme, lying between 750 and 4500 for different time points and lattice steps $\Delta x$. The $L^2$ errors and the orders of convergence, estimated as in the previous section by halving the space step, are presented in Table~\ref{table_BGRWconv_sat_noniter}.

\begin{table}[!ht]
\caption{Convergence of the non-iterative BGRW solver for coupled saturated\\ flow and Monod reactions.}
\label{table_BGRWconv_sat_noniter}
\begin{tabular}{ c c c c c c c  c c c c}
\hline
   & $\Delta x$ & $\|\psi-\psi_{m}\|$ & EOC & $\|c_1-c_{1m}\|$ & EOC  & $\|c_2-c_{2m}\|$ & EOC \\
  \hline
   &  2.00e-1  & 1.76e-01 & --   & 1.92e-02 & --   & 2.55e-02 & --   &\\
   &  1.00e-1  & 4.44e-02 & 1.99 & 4.83e-03 & 1.99 & 5.88e-03 & 2.11 &\\
   &  5.00e-2  & 1.21e-02 & 1.87 & 1.27e-03 & 1.93 & 1.47e-03 & 2.00 &\\
   &  2.50e-2  & 7.22e-03 & 0.75 & 5.57e-04 & 1.19 & 4.88e-04 & 1.59 &\\
  \hline
\end{tabular}
\end{table}

For comparison, the same problem is solved with the iterative BGRW $L$-scheme from Section~\ref{sec:iter}. For stabilization parameter $L=1$ and tolerances $\varepsilon_a=\varepsilon_r=10^{-6}$ in (\ref{eq:5}), the convergence of the liniarization procedure for different $\Delta x$ values is already established after about 25 iterations. The corresponding convergence results are shown in Table~\ref{table_BGRWconv_sat_iter}. We note that there is a similar convergence behavior of the iterative and non-iterative solvers. The results shown in Tables~\ref{table_BGRWconv_sat_noniter}~and~\ref{table_BGRWconv_sat_iter} can be further improved, at the price of larger computation time, by increasing the stabilization parameter $L$, which leads to smaller error norms and larger EOC \cite[Table 4]{Suciuetal2021}.

Further, we solve the Eqs.~(\ref{eq:Mobile})-(\ref{eq:Immobile}) for reactive transport with a given analytical solution of the flow problem. The flow velocity is obtained by differentiating the pressure head (\ref{eq:p_sat}), according to (\ref{eq:Q}). The results obtained with the non-iterative BGRW scheme, presented in Table~\ref{table_BGRWconv_sat_analytic}, show a much improved order of convergence. Hence, the smaller EOC values of the concentration solutions shown Tables~\ref{table_BGRWconv_sat_noniter}~and~\ref{table_BGRWconv_sat_iter} can be attributed to the coupling and the limited accuracy of the flow solver. We mention that almost the same convergence results (not shown here) were obtained with the iterative scheme but the computation was several times slower (30 s vs. 4 s for the non-iterative solver). Hence, the non-iterative solver is more appropriate to solve Monod-type problems in case of prescribed velocity.

\begin{table}[!ht]
\caption{Convergence of the iterative solver based on the BGRW $L$-scheme for\\ coupled saturated flow and Monod reactions.}
\label{table_BGRWconv_sat_iter}
\begin{tabular}{ c c c c c c c  c c c}
\hline
   & $\Delta x$ & $\|\psi-\psi_{m}\|$ & EOC & $\|c_1-c_{1m}\|$ & EOC  & $\|c_2-c_{2m}\|$ & EOC \\
  \hline
   &  2.00e-1  & 2.94e-01 & --   & 1.50e-02 & --   & 1.59e-02 & --   &\\
   &  1.00e-1  & 7.36e-02 & 2.00 & 3.86e-03 & 1.96 & 4.06e-03 & 1.96 &\\
   &  5.00e-2  & 1.90e-02 & 1.95 & 1.06e-03 & 1.86 & 1.11e-03 & 1.87 &\\
   &  2.50e-2  & 7.27e-03 & 1.39 & 4.27e-04 & 1.32 & 3.57e-04 & 1.63 &\\
  \hline
\end{tabular}
\end{table}

\begin{table}[!ht]
\caption{Convergence of the non-iterative BGRW scheme for given\\ analytical flow velocity and Monod reactions.}
\label{table_BGRWconv_sat_analytic}
\begin{tabular}{ c c c c c c c  c c c}
\hline
   & $\Delta x$ &  $\|c_1-c_{1m}\|$ & EOC  & $\|c_2-c_{2m}\|$ & EOC \\
  \hline
   &  2.00e-1  & 2.67e-02 & --   & 4.12e-02 & -- &\\
   &  1.00e-1  & 5.32e-03 & 2.32 & 7.38e-03 & 2.48 &\\
   &  5.00e-2  & 1.22e-03 & 2.13 & 1.56e-03 & 2.24 &\\
   &  2.50e-2  & 2.84e-04 & 2.10 & 3.49e-04 & 2.16 &\\
  \hline
\end{tabular}
\end{table}

\subsection{Saturated/unsaturated flows and Monod reactions}
\label{sec:unsat_flow}

To model both unsaturated and saturated flow conditions, we consider the pressure head solution
\begin{equation}
\psi_m(x,z,t) = -tx(2-x)z(3-z)-1+x/2+y/3.\label{eq:p_unsat}
\end{equation}
As an example of nontrivial parameters which cause the degeneracy of the Richards equation, we consider the variable water content and hydraulic conductivity given by
\begin{equation} \label{eq:theta}
\theta(\psi) = \begin{cases} 1/(3.3333-\psi) &\psi < 0 \\
0.3 &\psi \geq 0,
\end{cases}
\end{equation}
\begin{equation} \label{eq:K}
K(\Theta(\psi)) = \begin{cases} 0.05[1-(0.3-\theta(\psi))] &\psi < 0 \\
0.05 &\psi \geq 0.
\end{cases}
\end{equation}
The analytical solutions for the mobile species are again given by (\ref{eq:c1})-(\ref{eq:c2}) and the reactions are governed by the same Monod model as in Section~\ref{sec:sat_flow}.

For the beginning, we investigate the convergence behavior of the BGRW $L$-scheme in conditions of strictly unsaturated flow by removing the last three terms from (\ref{eq:p_unsat}), with variable $\theta$ given by (\ref{eq:theta}) and a constant hydraulic conductivity $K=0.05$. With the same parameters $L$, $\varepsilon_a$, and $\varepsilon_r$ as for the saturated case, the iterative procedure already converges after 50 to 100 iterations for the flow solution and 23 to 28 iterations for the concentration solutions. The convergence behavior shown in Table~\ref{table_BGRWconv_unsat} is similar to that in the saturated case shown in Table~\ref{table_BGRWconv_sat_iter}. The noticeable difference from the saturated case is given by the much shorter computing time of cca. 5 min, vs. cca 3 h in saturated case. Since, as we have seen above, the computation of the reactive transport alone takes only a few seconds, this difference is caused by the slower convergence of the flow solver in the saturated case.

\begin{remark}
\label{rem:unsat-noniter}
The non-iterative BGRW scheme from Section~\ref{sec:noniter} can be adapted for variable water content $\theta$ by adding to the right-hand side of Eq.~(\ref{eq:BGRW1}) the source term
$f_{i,j,k}=-c_{i,j,k}(\theta_{i,j,k+1}-\theta_{i,j,k})/\theta_{i,j,k}$. Performing a numerical test similar to that presented in Table~\ref{table_BGRWconv_unsat} leads to relative errors of the numerical solutions for the mobile concentrations $c_1$ and $c_2$ of about $15\%$, even if the spatial step $\Delta x$ is further decreased from $2.50\cdot 10^{-2}$ to $3.125\cdot 10^{-3}$. For comparison, the relative errors of the solutions provided by the non-iterative scheme for $\Delta x=2.50\cdot 10^{-2}$ in the convergence test for saturated flow presented in Table~\ref{table_BGRWconv_sat_noniter} are about $0.08\%$ for $c_1$ and $0.05\%$ for $c_2$. Therefore, while the non-iterative scheme is accurate and several times faster than the iterative BGRW scheme in case of saturated flows, it is not appropriate for unsaturated flow conditions.
\end{remark}

\begin{table}[!ht]
\caption{Convergence of the iterative solver based on the BGRW $L$-scheme for\\ coupled unsaturated flow and Monod reactions.}
\label{table_BGRWconv_unsat}
\begin{tabular}{ c c c c c c c  c c c}
\hline
   & $\Delta x$ & $\|\psi-\psi_{m}\|$ & EOC & $\|c_1-c_{1m}\|$ & EOC  & $\|c_2-c_{2m}\|$ & EOC \\
  \hline
   &  2.00e-1  & 2.94e-01 & --   & 2.37e-02 & --   & 3.54e-02 & --   &\\
   &  1.00e-1  & 7.40e-02 & 1.99 & 5.94e-03 & 2.00 & 9.36e-03 & 1.92 &\\
   &  5.00e-2  & 1.92e-02 & 1.95 & 1.67e-03 & 1.83 & 2.64e-03 & 1.83 &\\
   &  2.50e-2  & 6.91e-03 & 1.47 & 5.56e-04 & 1.59 & 7.57e-04 & 1.80 &\\
  \hline
\end{tabular}
\end{table}

\begin{table}[!ht]
\caption{Convergence of the iterative flow solver from \cite{Suciuetal2021} for solutions of degenerate\\ Richards equation describing transitions from unsaturated to saturated flow regimes.}
\label{table_BGRWconv_flow}
\begin{tabular}{ c c c c c c c  c c c}
\hline
   & $\Delta x$ & $\|\psi-\psi_{m}\|$ & EOC & $\|U-U_{m}\|$ & EOC  & $\|V-V_{m}\|$ & EOC \\
  \hline
   &  2.00e-1  & 6.62e-01 & --   & 1.85e-02 & --   & 9.17e-03 & --   &\\
   &  1.00e-1  & 1.69e-01 & 1.97 & 6.87e-03 & 1.43 & 2.86e-03 & 1.68 &\\
   &  5.00e-2  & 4.22e-02 & 2.00 & 2.40e-03 & 1.52 & 9.36e-04 & 1.61 &\\
   &  2.50e-2  & 1.06e-02 & 2.00 & 8.15e-04 & 1.56 & 3.06e-04 & 1.61 &\\
  \hline
\end{tabular}
\end{table}

\begin{table}[!ht]
\caption{Convergence of the concentration solutions of the\\ BGRW $L$-scheme coupled with the flow solver from \cite{Suciuetal2021}\\ for saturated/unsaturated flow and Monod reactions.}
\label{table_BGRWconv_sat_unsat}
\begin{tabular}{ c c c c c c c  c c c}
\hline
   & $\Delta x$ & $\|c_1-c_{1m}\|$ & EOC  & $\|c_2-c_{2m}\|$ & EOC \\
  \hline
   &  2.00e-1  & 3.14e-02 & --   & 5.84e-02 & --   &\\
   &  1.00e-1  & 1.02e-02 & 1.62 & 1.88e-02 & 1.64 &\\
   &  5.00e-2  & 2.86e-03 & 1.83 & 5.23e-03 & 1.84 &\\
   &  2.50e-2  & 7.32e-04 & 1.97 & 1.34e-03 & 1.97 &\\
  \hline
\end{tabular}
\end{table}

In the following, we consider the challenging situation where the hydraulic conductivity $K$ is given by the step function (\ref{eq:K}). Now, both $\theta$ and $K$ take constant values in saturated regions and Richards equation degenerates into an elliptic equation. The one-way coupling of the two-dimensional flow and reactive transport $L$-schemes is tested by comparisons with the exact concentration solutions (\ref{eq:c1})-(\ref{eq:c2}) and the flow solution (\ref{eq:p_unsat}), which captures the transition from unsaturated to saturated regime. The setup of the reaction system is the same as that used in the previous sections to test particular cases of the general problem. The system of coupled equations (\ref{eq:Q})-(\ref{eq:Immobile}), with source terms in Eqs.~(\ref{eq:Richards})~and~(\ref{eq:Mobile}) derived from the analytical solution, is solved numerically with the flow and transport $L$-schemes. With the same $\varepsilon_a$, and $\varepsilon_r$ as used above and stabilization parameter set to $L=10$ in the flow solver and to $L=1$ in the transport solver, the desired accuracy of the iteration procedure is attained after $10^3$ to $2\cdot 10^4$ iterations of the flow scheme and 30 to 65 iterations of the transport scheme. The total computation time for the four discretization levels considered in this test is about 6 hours.

Convergence results for the flow and transport solvers are presented in Tables~\ref{table_BGRWconv_flow}~and~\ref{table_BGRWconv_sat_unsat}, respectively. One remarks a tendency towards second order of convergence of the pressure head and the concentration solutions. The $L^2$ error norms of the velocity solutions are smaller by one order of magnitude than those for pressure and concentration. However, their convergence order is smaller, in the range from 1.4 to 1.7.

We note that the influence of the stabilization parameter $L$ on accuracy and order of convergence estimations (cf. \cite[Sect. 5.1]{Suciuetal2021}) has not been investigated in the present study. Also, we did not tested the unbiased GRW non-iterative scheme presented in Section~\ref{sec:ugrw} by comparisons with analytical solutions of Monod problems. As seen in Section~\ref{sec:const_flow}, the implementation of the boundary conditions in unbiased GRW schemes may distort the estimated order of convergence. The performance of the GRW scheme will be evaluated in Section~\ref{sec:bioSat} below through comparisons with BGRW solutions.

\section{Biodegradation in saturated flows}
\label{sec:bioSat}

In the following we consider a saturated regime with completely decoupled flow and reactive transport processes. The water content $\theta$ equals the porosity and enters as a constant in the BGRW and GRW parameters given by relations (\ref{eq:ParBGRW1}) and (\ref{eq:GRW2}), respectively. The two GRW algorithms will be used to build up numerical solutions to biodegradation problems in conditions of advection-dominated transport and heterogeneous  structure of the aquifer system.

\subsection{Biodegradation and advection-dominated transport}
\label{sec:bioAdv}

The artificial diffusion induced by numerical schemes is a serious challenge in simulations of biodegradation processes. We refer to \cite{Raduetal2011} for a comprehensive study on numerical diffusion for different discretization schemes. An artificial mixing of the species leads to an overprediction of the degradation or even to a complete degradation of the contaminant, contrary to what would be observed in practice. Mixed finite element methods and higher order schemes, together with stabilization techniques were proposed in the past to cope with the advection-dominated transport. Some drawbacks of these approaches have been also reported, such as the occurrence of negative concentrations (in mixed finite element approaches \cite{Brunneretal2012}) or troubles with the construction of efficient iterative solvers (when using quadratic finite element schemes and streamline upwind stabilization \cite{BauseandKnabner2004}).

As an alternative, we propose GRW solutions to the benchmark problem used in the papers cited above. Besides being free of numerical diffusion (provided that $\text{P\'{e}}\le 2$ in case of BGRW algorithm), the GRW solvers working with indivisible particles do not produce negative concentrations, the conversion between particle densities and concentrations being straightforward, by using the same number of particles as molecules involved in reactions. Moreover, for saturated flows conditions, as in case of this benchmark problem, we have the advantage of using non-iterative BGRW/GRW schemes which, as seen in Section~\ref{sec:sat_flow} above, are faster than the iterative schemes.

The problem is formulated in the domain $\Omega=(0,0.5)\times(0,1)$. The reactive transport of the mobile species is governed by Eqs.~(\ref{eq:Mobile}) with parameters specified by $\theta=0.2$, $U=0$, $V=-1$, $D_1=D_2=D=10^{-4}$. The Monod reaction term $\mu$ is modeled by Eq.~(\ref{eq:Monod}) with $M_1=M_2=0.1$, the biomass is kept constant, $c_3=1$, and the reaction terms in Eqs.~(\ref{eq:Mobile}) are given by $R_1=\theta\alpha_1\mu$ $R_2=\theta\alpha_2\mu$, with $\alpha_1=5$ and $\alpha_2=0.5$. The final simulation time is $T=5$, when one expects that the process reaches the steady state \cite{BauseandKnabner2004}. As initial conditions, one chooses $c_1(x,z,0)=0$ and $c_2(x,z,0)=0.1$ in $\Omega$. The contaminant (electron donor) is injected at the middle part of the inflow boundary, where no electron acceptor (oxygen) is present, i.e., $c_1(x,z,t)=1$ and $c_2(x,z,t)=0$
for $(x,z,t)\in(0.225,0.275)\times\{1\}\times(0,T]$. Further, one sets $c_1(x,z,t)=0$ and $c_2(x,z,t)=0.1$ for $(x,z,t)\in(0,0.5)\backslash(0.225,0.275)\times\{1\}\times(0, T]$. Finally, the condition of vanishing concentration gradients is imposed on the left, right, and bottom boundaries.

\begin{figure}
\begin{minipage}[t]{0.45\linewidth}\centering
\includegraphics[width=\linewidth]{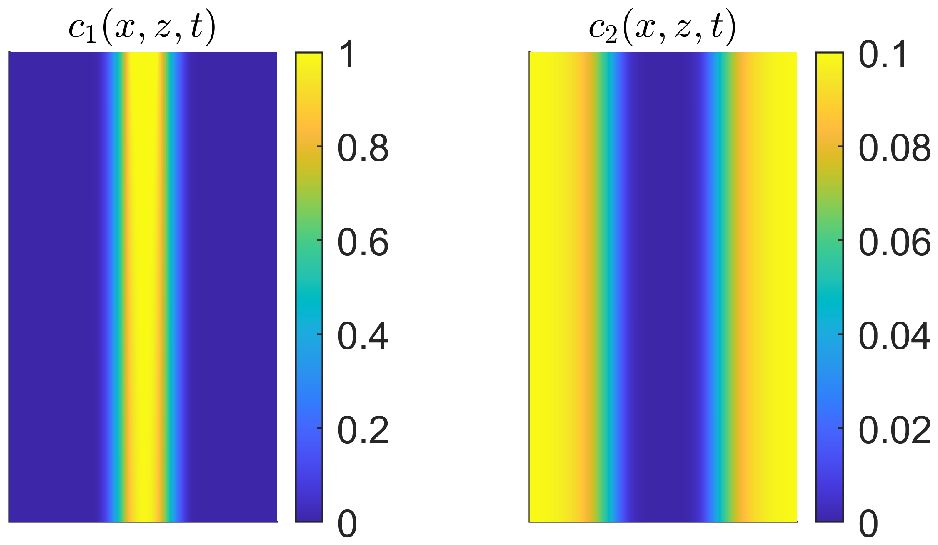}
\caption{\label{fig_1_numdiff}Concentration profiles at the final time $T=5$ for $D=10^{-4}$.}
\end{minipage}
\hspace*{0.1in}
\begin{minipage}[t]{0.45\linewidth}\centering
\includegraphics[width=\linewidth]{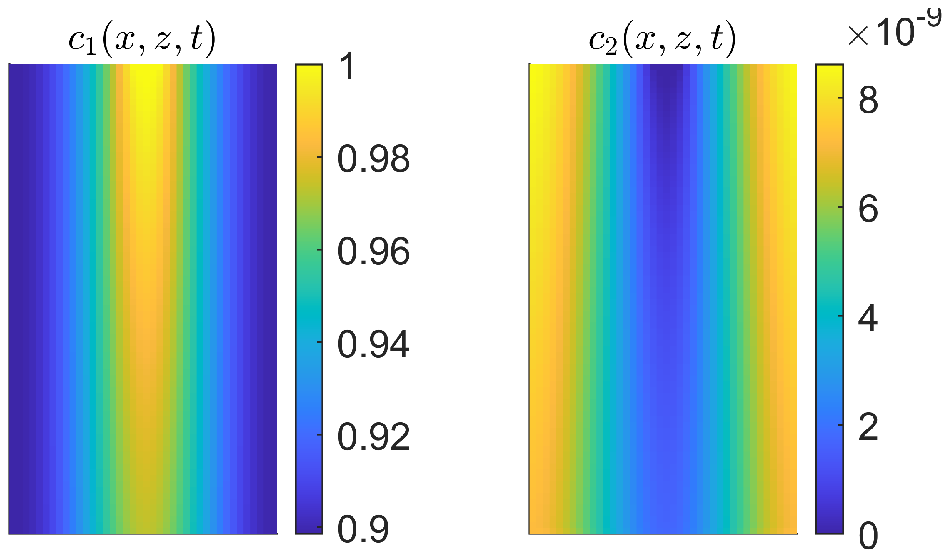}
\caption{\label{fig_2_numdiff}Concentration profiles at the final time $T=5$ for $D=2.5\cdot 10^{-3}$.}
\end{minipage}
\end{figure}

The BGRW computations for the present setup require an extremely fine discretization ($\Delta x=2\cdot 10^{-4}$ and $\Delta t=10^{-5}$ for $\text{P\'{e}}=2$) which would render the computation until the final time $T=5$ prohibitive. Therefore, we first restrict the final time to $T=1$ and compare the BGRW solution to that provided by the unbiased GRW solver. The latter is obtained with $\Delta x=2\cdot 10^{-3}$ ($\text{P\'{e}}=20$), $\Delta t=4\cdot 10^{-4}$, and GRW parameters $v=1$, $d=1$. The comparison is done in terms of concentrations averaged over lattice sites, $\langle c_{\nu}\rangle,\; \nu=1,2$, by relative errors $\varepsilon_{c_{\nu}}=(\langle c_{\nu}\rangle-\langle c_{\nu}^{b}\rangle)/\langle c_{\nu}^{b}\rangle$, where $c_{\nu}^{b}$ is the BGRW solution. The errors obtained, $\varepsilon_{c_1}=-0.0125$ and $\varepsilon_{c_2}=-0.0016$, validate the GRW solver for the parameters of the benchmark problem described above and $T=1$. Further, the problem is solved for $T=5$ with the unbiased GRW solver. The electron donor and acceptor concentrations at the final time shown in Fig.~\ref{fig_1_numdiff} are similar to those obtained by quadratic finite element in \cite[Fig. 2]{BauseandKnabner2004}. One remarks that the contaminant is advected by the constant flow, with negligible diffusion and biodegradation.

Even though the GRW solver does not produce numerical diffusion, the effect of an artificial mixing which causes an overprediction of the biodegradation can be illustrated by increasing the diffusion coefficient. To simulate an enhancement of the mixing process, we repeat the computations after setting the diffusion coefficient to $D=2.5\cdot 10^{-3}$. The results are presented in Fig.~\ref{fig_2_numdiff}. The concentration profiles are now much more diffusive, the oxygen is largely consumed, as shown by the small values $c_2$, and the more intense biodegradation in the lower part of the domain leads to the underestimation of the contaminant concentration $c_1$ at the outflow boundary.

\subsection{Influence of the heterogeneity on the biodegradation process}
\label{sec:bioHeterogAdv}

As seen in Section~\ref{sec:bioAdv} above, in case of localized reaction zones the increase of the diffusion coefficient results in an increased biodegradation of the electron donor contaminant and consumption of the oxygen. The impact of the artificial increase of the diffusion in numerical schemes can be severe especially in the case of biodegradation controlled by transverse mixing in heterogeneous media \cite{Cirpkaetal1999}. Classical approaches to cope with the numerical diffusion in advection-dominated problems are based on streamtube reactive transport models \cite{Cirpkaetal1999}, local mesh refinement \cite{Klofkornetal2002}, or higher-order finite element used together with stabilization techniques \cite{BauseandKnabner2004}. In the following we propose GRW solutions to transport and biodegradation problems similar to those considered in these papers and investigate the influence of varying the parameters describing the heterogeneity.

We consider a problem similar to that formulated in \cite{Cirpkaetal1999}, with a two-dimensional domain $\Omega=(0,80)\times(0,20)$, a uniform initial distribution of biomass, a continuous injection of contaminant, constant concentration of oxygen in the water entering over the inflow boundary and in the domain, excepting the injection well. The water flow is driven by the gradient of the piezometric pressure and a variable hydraulic conductivity modeled as a realization of a random field. In the following, the distances are measured in meters, the time in days, and the concentration in moles.

To facilitate the computation of Monte Carlo ensembles of solutions, the flow velocity is approximated by a Kraichnan procedure \cite[Appendix C.3.2.2]{Suciu2019} as sum of 100 cosine random modes, with amplitudes,  wavenumbers and phased determined by the parameters of the random hydraulic conductivity $K$. The log-hydraulic conductivity field $\ln K$ is a statistically homogeneous random function with isotropic Gaussian correlation. Its variance belongs to the set $\sigma^2\in\{0.1, 0.5,1\}$ and for its correlation length we consider $\lambda=1$ and $\lambda=2$. The mean conductivity is set to $\langle K\rangle=25$. With a decrease of $0.8$ of the pressure head between the inflow and outflow boundaries and a gradient $J=-0.8/80$, the mean velocity becomes $U=-J\langle K\rangle=0.25$. The
dispersion is parameterized by a constant coefficient $D=0.01$. Ensemble averages are estimated over 256 realizations of the $\ln K$ field. These parameters provide fairly good approximations of the random Darcy velocity field and of the stochastic advection-dispersion process \cite{Schwarzeetal2001,Eberhardetal2007}. Since the porosity does not play a major role in saturated flows it is set to $\theta=1$.

The setup of the chemical reactions closely follows that used in \cite{Cirpkaetal1999,BauseandKnabner2004}, with the only difference that, instead of mg/l, the concentrations are now given in moles. The parameters of the Monod model are given by $M_1=2,\; M_2=0.2,\; Y=0.09,\; K_d=0.05$, $\mu_{max}=5$, $\alpha_1=1,\; \alpha_2=3$ \cite{Cirpkaetal1999}. Additionally, following \cite{BauseandKnabner2004} we consider a growth limitation in the Monod model with $c_{3_{max}}=1$ and $\gamma=1$ in Eq.~(\ref{eq:Immobile}). The initial concentration of the contaminant is set to $c_1=2$ at lattice sites inside a square of side equal to 2 centered at $x=2.5$ and $y=10$ and to $c_1=0$ otherwise. The oxygen concentration is initially set to $c_2=5$ in $\Omega$, excepting the injection region where it is $c_2=0$. The continuous injection is simulated, as in \cite{BauseandKnabner2004}, by resetting after every time step both $c_1$ and $c_2$ to their initial values in the injection region. The initial concentration of the biomass is set to $c_3=0.001$ in $\Omega$. In \cite{BauseandKnabner2004}, the boundary conditions for the mobile species are set to the initial values of $c_1$ and $c_2$ on the top, left and bottom boundaries. To avoid boundary effects in the unbiased GRW solver (see Section~\ref{sec:const_flow}) which could lead to vanishing concentration $c_2$ close to boundaries we also fix $c_2$ to the initial value on a stripe of width 1.5 along the same boundaries. No boundary conditions are needed on the right boundary, because if no particles are introduced from from exterior of $\Omega$ both the GRW and BGRW algorithms simulate the physical outflow situation by removing the particles crossing the boundary.

\begin{figure}[hbt!]
\begin{minipage}[t]{0.45\linewidth}\centering
\includegraphics[width=\linewidth]{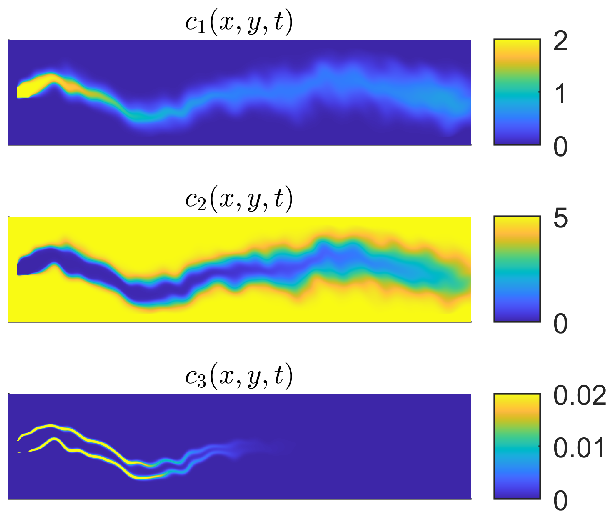}
\caption{\label{fig_bgrw}Concentration profiles of electron donor ($c_1$), electron acceptor ($c_2$), and biomass ($c_3$) at the final time $T$, obtained with the BGRW solver for $\lambda=2$, $\sigma^2=0.5$.}
\end{minipage}
\hspace*{0.1in}
\begin{minipage}[t]{0.45\linewidth}\centering
\includegraphics[width=\linewidth]{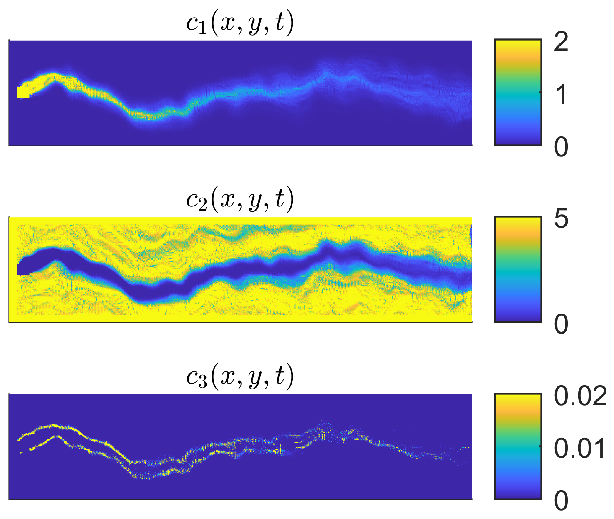}
\caption{\label{fig_grw}Concentration profiles of electron donor ($c_1$), electron acceptor ($c_2$), and biomass ($c_3$) at the final time $T$, obtained with the GRW solver for $\lambda=2$, $\sigma^2=0.5$.}
\end{minipage}
\end{figure}

In Fig.~\ref{fig_bgrw} we present the concentration profiles of the three species at the final time $T=300$ days computed with the non-iterative BGRW solver for a fixed realization of the velocity field with $\lambda=2$ and $\sigma^2=0.5$. The results are obtained on a lattice with $\Delta x=\Delta y=0.0333$ which ensures a maximum P\'{e}clet number P\'{e}=1.75 verifying the relations (\ref{eq:RestrBGRW}). The corresponding time step, according to (\ref{eq:ParBGRW1}) is $\delta t=0.0231$. The BGRW solution is further used as reference for comparison with that provided by the non-iterative unbiased GRW solver for the same realization of the velocity field. The latter is obtained by using a much coarser discretization,  with $\Delta x=\Delta y=0.1250$ (maximum P\'{e}=6.57) and $\delta t=1$ corresponding to $d=2$ and $u=U\delta t/\delta x=2$ in (\ref{eq:GRW2}). With these, the GRW solver achieves a speedup of $10^4$ over the BGRW solver.

The GRW solutions shown in Fig.~\ref{fig_grw} reproduce the general trend of the concentration profiles of the BGRW solution (Fig.~\ref{fig_bgrw}) and the shape of the reaction front. However, as indicated by the middle panel of Fig.~\ref{fig_grw}, there is an overshoot of the electron-donor particles in the GRW solution, which produces reactions beyond the reaction fronts predicted by the reference BGRW solution shown in Fig.~\ref{fig_bgrw}. The more extended profiles of the biomass concentration $c_3$ shown in the lower panel of Fig.~\ref{fig_grw} also can be explained by the numerical overshooting artifact. The effect of the overshooting can be quantitatively significant. In the present case, the relative errors of the lattice-averaged GRW solutions with respect to the BGRW solution are $\varepsilon_{c_{1}}=-0.3109$, $\varepsilon_{c_{2}}=0.0520$, and $\varepsilon_{c_{3}}=0.8256$, that is, one order of magnitude larger than in case of the uniform velocity field analyzed in Section~\ref{sec:bioAdv}. The increased errors can be attributed to the spatial variation of the advection velocity, requiring larger parameters $d$ and $v$ which give rise to enhanced overshooting. The negative value $\varepsilon_{c_{1}}$ quantifies the overprediction of the degradation caused by overshooting.

\begin{figure}
\begin{minipage}[t]{0.45\linewidth}\centering
\includegraphics[width=\linewidth]{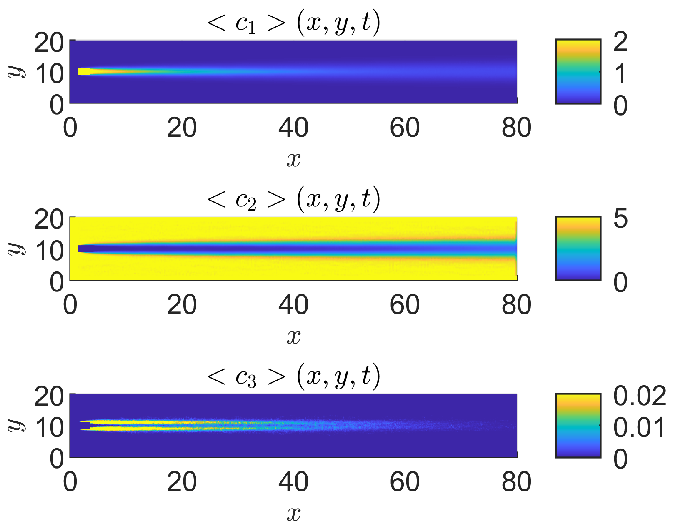}
\caption{\label{fig_meanGRW}Ensemble average concentrations at the final time $T$ for $\lambda=1$, $\sigma^2=0.1$.}
\end{minipage}
\hspace*{0.1in}
\begin{minipage}[t]{0.45\linewidth}\centering
\includegraphics[width=\linewidth]{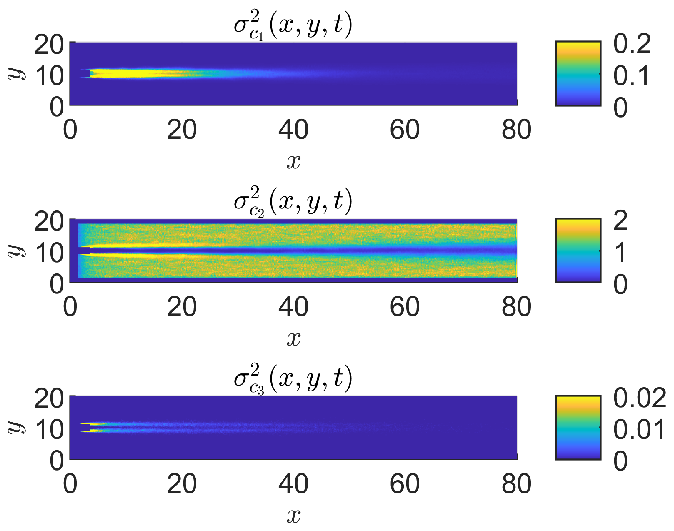}
\caption{\label{fig_varGRW}Concentration variances at the final time $T$ for $\lambda=1$, $\sigma^2=0.1$.}
\end{minipage}
\end{figure}

Keeping in mind that the GRW solutions overpredict the biodegradation, they can still be used to investigate the effect of heterogeneity through Monte Carlo simulations. The species concentrations averaged over 256 realizations of the $\ln K$ field with $\lambda=1$, $\sigma^2=0.1$ and the corresponding variances are shown in Figs.~\ref{fig_meanGRW}~and~\ref{fig_varGRW}. The averaged concentrations shown in Fig.~\ref{fig_meanGRW} are similar to the results obtained with uniform velocity. But the relative errors of the lattice-averaged mean concentrations from Fig.~\ref{fig_meanGRW} with respect to the BGRW solutions for a constant longitudinal velocity set to $U=0.25$ are close to $\varepsilon_{c_{\nu}}$ in the comparison of the GRW/BGRW solutions for a single realization of the $\ln K$ field presented above. The large variance $\sigma^{2}_{c_2}$ in Fig.~\ref{fig_varGRW}, consistent with the inhomogeneity of the electron acceptor concentration $c_2$ in a single realization shown in Fig.~\ref{fig_grw}, is also partially due the overshooting in the GRW solver. Averages over lattice sites of the ensemble-mean species concentrations and variances for different correlation lengths $\lambda$ and variances $\sigma^2$ of the $\ln K$ field, computed at the final time $T$, are presented in Table~\ref{table_MC}. The mean species concentrations show some variation, of about 22\% ($\langle c_1\rangle$), 25\% ($\langle c_2\rangle$), and 12\% ($\langle c_3\rangle$), but there is no systematic dependence on $\lambda$ and $\sigma^2$. The variances $\sigma^{2}_{c_{\nu}}$ slightly increase with $\sigma^2$ but do not show a systematic dependence on $\lambda$.

\begin{table}[!ht]
\caption{Monte Carlo results. Ensemble mean species  concentrations\\ and variances averaged over lattice sites.}
\label{table_MC}
\begin{tabular}{ c c c c c c c c c c c c}
\hline
   & $\lambda$ & $\sigma^2$ & $<c_1>$ & $<c_2>$ & $<c_3>$ & $\sigma^{2}_{c_1}$ & $\sigma^{2}_{c_2}$ & $\sigma^{2}_{c_3}$\\
  \hline
   &  1  & 0.1  & 9.07e-2 & 4.08e+0 & 1.57e-3 & 1.03e-2 & 1.03e+0 & 2.60e-4 \\
   &  1  & 0.5  & 8.66e-2 & 3.04e+0 & 1.69e-3 & 2.38e-2 & 2.30e+0 & 3.90e-4 \\
   &  1  & 1.0  & 7.76e-2 & 4.00e+0 & 1.76e-3 & 3.07e-2 & 4.27e+0 & 4.93e-4 \\
   &  2  & 0.1  & 9.55e-2 & 4.08e+0 & 1.57e-3 & 1.93e-2 & 1.36e+0 & 2.47e-4 \\
   &  2  & 0.5  & 1.00e-1 & 4.02e+0 & 1.68e-3 & 4.09e-2 & 2.33e+0 & 3.50e-4 \\
   &  2  & 1.0  & 9.98e-2 & 4.00e+0 & 1.79e-3 & 5.17e-2 & 3.27e+0 & 4.25e-4 \\
  \hline
\end{tabular}
\end{table}

\begin{figure}[p]
\centering
\includegraphics[width=\linewidth]{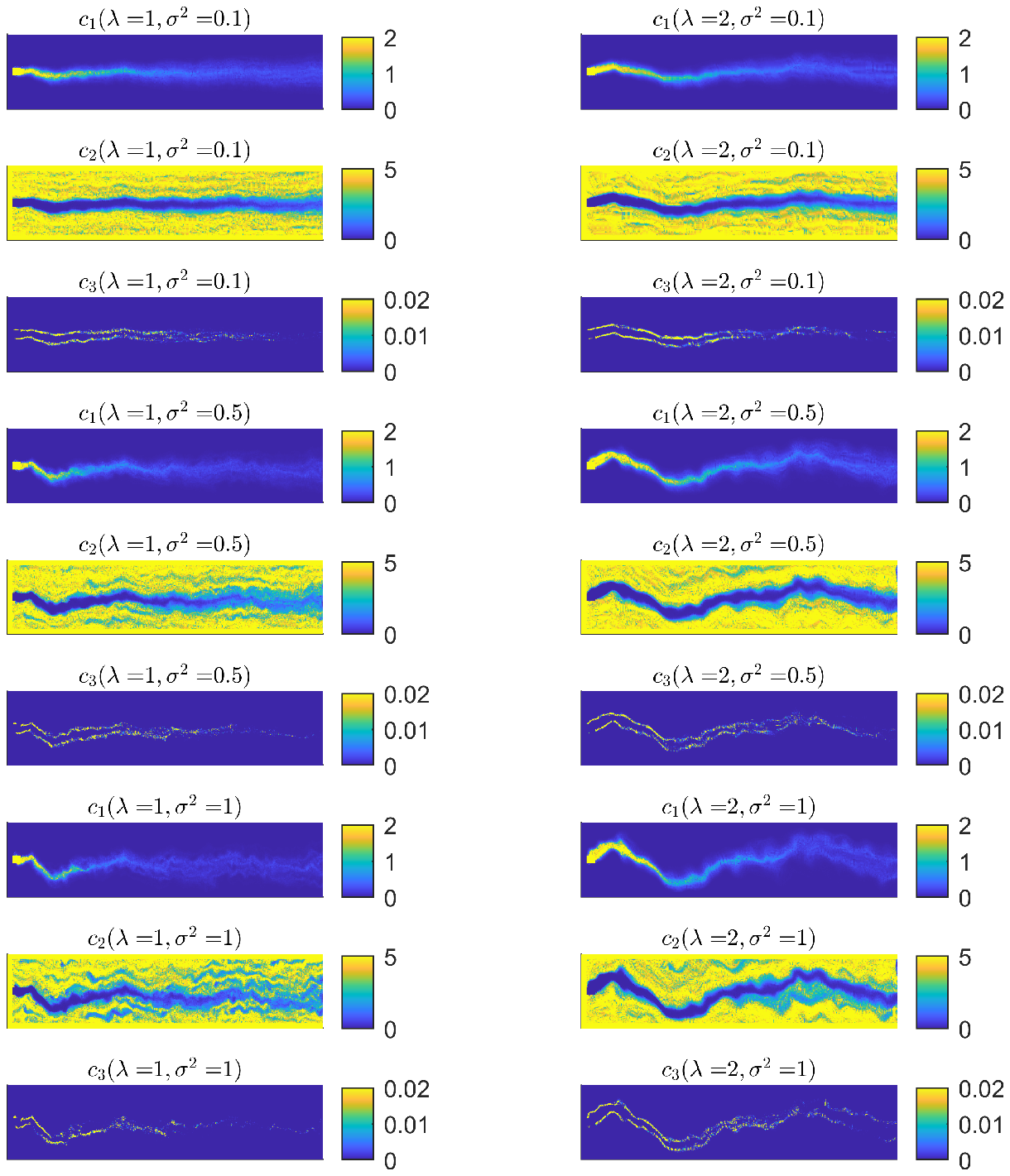}
\caption{\label{fig_comp_grw}Concentration profiles of electron donor ($c_1$), electron acceptor ($c_2$), and biomass ($c_3$) at the final time $T$, obtained with the unbiased GRW solver for $\lambda\in\{1, 2\}$, $\sigma^2\in\{0.1, 0.5, 1\}$.}
\end{figure}

The comparison of $c_1$, $c_2$, and $c_3$ for the first realization in the Monte Carlo ensembles computed for the six combinations of $\lambda$ and $\sigma^2$, presented in Fig.~\ref{fig_comp_grw}, shows that the spatial variability of the velocity field mainly influences the meandering of the reaction fronts and the amount of overshooting. As seen in both columns of panels in Fig.~\ref{fig_comp_grw}, for fixed $\lambda$ the meandering is more pronounced as the variance $\sigma^2$ increases. Comparing the two columns, we remark an increase in meandering amplitude as the correlation length $\lambda$ of the $\ln K$ filed increases. We also remark that the amount of overshooting indicated by zones of low concentrations $c_2$ of the electron acceptor outside the reaction fronts, increases with both $\lambda$ and $\sigma^2$.

\begin{remark}\label{rem:overshoot}
We note that the overshooting is an issue of concern for both unbiased GRW and particle tracking approaches. The GRW rule (\ref{eq:GRW1}) moves particles in space, between the lattice sites, by summing discrete advective displacements and diffusive jumps, e.g., $v\pm d$, in the same way as a numerical solution of the It\^o equation, which is the mathematical model of the PT procedure. The only differences are, first, that $v$ and $d$ are discrete quantities which measure displacements between lattice sites, and second, instead of computing ensembles of PT trajectories, the GRW scheme counts the number of particles at lattice sites which approximate the solution of the partial differential equation to be solved. In both schemes, displacements of computational particles over space points where the problem data may have significant variations, e.g., diffusion fronts, are hardly avoidable. For instance, in the semi-Lagrangian method proposed by \cite{Natarajanetal2021}, to prevent a tracked particle from leaving the computational subdomain, the time step is restricted such that the sum of advective and dispersive displacements is smaller than the minimum grid spacing. Further, the polynomial solution of the transport problem associated to the particle has to be interpolated to the nodes of the Chebyshev quadrature where a coupled spectral element method provides the advection velocity. For comparison, in the BGRW method, which also restricts the time step, the smooth solution is represented by a sufficiently large number of particles which jump only to first-neighbor lattice sites and the velocity is given at the same lattice sites (by another GRW solver in coupled flow and transport problems). One avoids in this way both the overshooting and the interpolation procedure.
\end{remark}

\section{Biodegradation in variably saturated soils}
\label{sec:deg_Rich}

In designing bioremediation strategies to clean up contaminated soils, the interplay between heterogeneity, variable saturation, and soil properties has to be taken into account in numerical models. As an illustration, we simulate in the following the infiltration and biodegradation of benzene in conditions of partial saturation for two types of soils.

We consider the same two-dimensional domain $\Omega=(0,2)\times(0,3)$ used in verification tests presented in Section~\ref{sec:verif}. The infiltration is driven by a variable pressure, given by $\psi(x,z,t\leq t_1)=-0.2+2.7t/t_1$ and $\psi(x,z,t> t_1)=0.2$, on the open part $\Gamma_1=(0.5,1.5)\times\{3\}$ of the top boundary and by communication with the fresh water at hydrostatic equilibrium on the open parts of the domain $\Gamma_{21}=\{0\}\times(0,0.5)$ and $\Gamma_{22}=\{2\}\times(0,0.5)$. The zero flow condition is imposed on the closed boundary $\Gamma\backslash(\Gamma_1\bigcup\Gamma_{21}\bigcup\Gamma_{22})$ and the initial pressure corresponds to hydrostatic equilibrium, $\psi(x,z,0)=0.5-z$ \cite{Radu2004}. Hereafter, the time is measured in days and the distances in meters. The intermediate time for the variable pressure condition is set to $t_1=1$ and the simulations are conducted for increasing times $T$.

The relationships defining the water content $\theta(\psi)$ and the hydraulic conductivity $K(\theta(\psi))$ are given by the van Genuchten-Mualem model
\begin{equation} \label{theta}
\Theta(\psi) = \begin{cases} \left(1+(-\alpha \psi)^n\right)^{-m} &\psi < 0 \\
1 &\psi \geq 0,
\end{cases}
\end{equation}
\begin{equation} \label{K}
K(\Theta(\psi)) = \begin{cases} K_{sat} \Theta(\psi)^{\frac{1}{2}} \left[1-\left(1-\Theta(\psi)^\frac{1}{m}\right)^m \right]^2 &\psi < 0 \\
K_{sat} &\psi \geq 0,
\end{cases}
\end{equation}
where $\Theta = (\theta - \theta_{res})/(\theta_{sat} - \theta_{res})$ is the normalized water content, $\theta_{res}$ and $\theta_{sat}$ are the residual and the saturated water content, $K_{sat}$ is the hydraulic conductivity of the saturated soil, and $\alpha$, $n$, $m=1-1/n$ are model parameters depending on the soil type.

We consider here two sets of soil parameters, presented in Table \ref{table:parametersRichards}. They correspond to a silt loam and a Beit Netofa clay, respectively, previously used in GRW simulations of coupled flow and nonreactive transport \cite[Sect. 5]{Suciuetal2021}.
\begin{table}[!ht] \centering
 \caption{Soil parameters}
 \label{table:parametersRichards}
 \begin{tabular}{ c c c }
  \hline
  & Silt loam & Beit Netofa clay \\
$\theta_{sat}$ & 0.396 & 0.446\\
$\theta_{res}$ & 0.131 & 0\\
$\alpha$ & 0.423 & 0.152\\
$n$ & 2.06 & 1.17\\
$K_{sat}$ & $4.96 \cdot 10^{-2}$ & $8.2\cdot 10^{-4}$\\
\hline
 \end{tabular}
\end{table}

The soil heterogeneity is modeled by a fixed realization of a random saturated hydraulic conductivity with mean $K_{sat}$. The latter is computed by the Kraichnan method described in \cite[Appendix C.3.1.2]{Suciu2019}, with an anisotropic Gaussian model of correlation lengths $\lambda_x=0.1$ $\lambda_z=0.01$, and variance $\sigma^2=4$. The dispersion is parameterized by the constant coefficient $D=0.003$.

At the initial time, the concentration of benzene (electron donor) is set to $c_1=2$ mole on the injection boundary $\Gamma_1$ and to $c_1=0$ otherwise, the concentration of oxygen (electron acceptor) is $c_2=0$ on $\Gamma_1$ and $c_2=5$ mole otherwise, and the biomass concentration is $c_3=0.001$ mole in $\Omega$. At $t\geq 0$, the concentrations $c_1$ and $c_2$ are kept at the initial values on the open boundaries and no flux conditions are imposed on $\Gamma\backslash(\Gamma_1\bigcup\Gamma_{21}\bigcup\Gamma_{22})$. With these, one simulates a continuous injection of benzene through $\Gamma_1$ and one assumes that the contaminated incoming water does not contain oxygen and the outgoing water is cleaned instantly at the contact with the fresh water on $\Gamma_{21}\bigcup\Gamma_{22}$. The reactions are governed by the Monod model with the same parameters as those used in Section~\ref{sec:bioHeterogAdv} for the simulations of reactive transport in aquifers.

The simulations are performed with the iterative BGRW $L$-scheme coupled with the GRW flow solver from \cite{Suciuetal2021}, by using the stabilization parameter $L=1$ in both the flow solver and the transport solvers, $\Delta x=\Delta z=0.05$, $\Delta t=8.68\cdot 10^{-5}$, for loam soil parameters, and $\Delta t=5.25\cdot 10^{-3}$ for clay soil parameters. The convergence of the $L$-schemes is achieved after tens to hundreds of iterations for $T<3$, with more iterations in loam soil simulations, and decreases to only a few iterations for larger $T$. As $T$ increases from $T=3$ to $T=120$, the computing time ranges between 8 min and 2.5 h for the loam soil and between 50 s and 8 min for the clay soil.

\begin{figure}[p]\centering
\includegraphics[width=\linewidth]{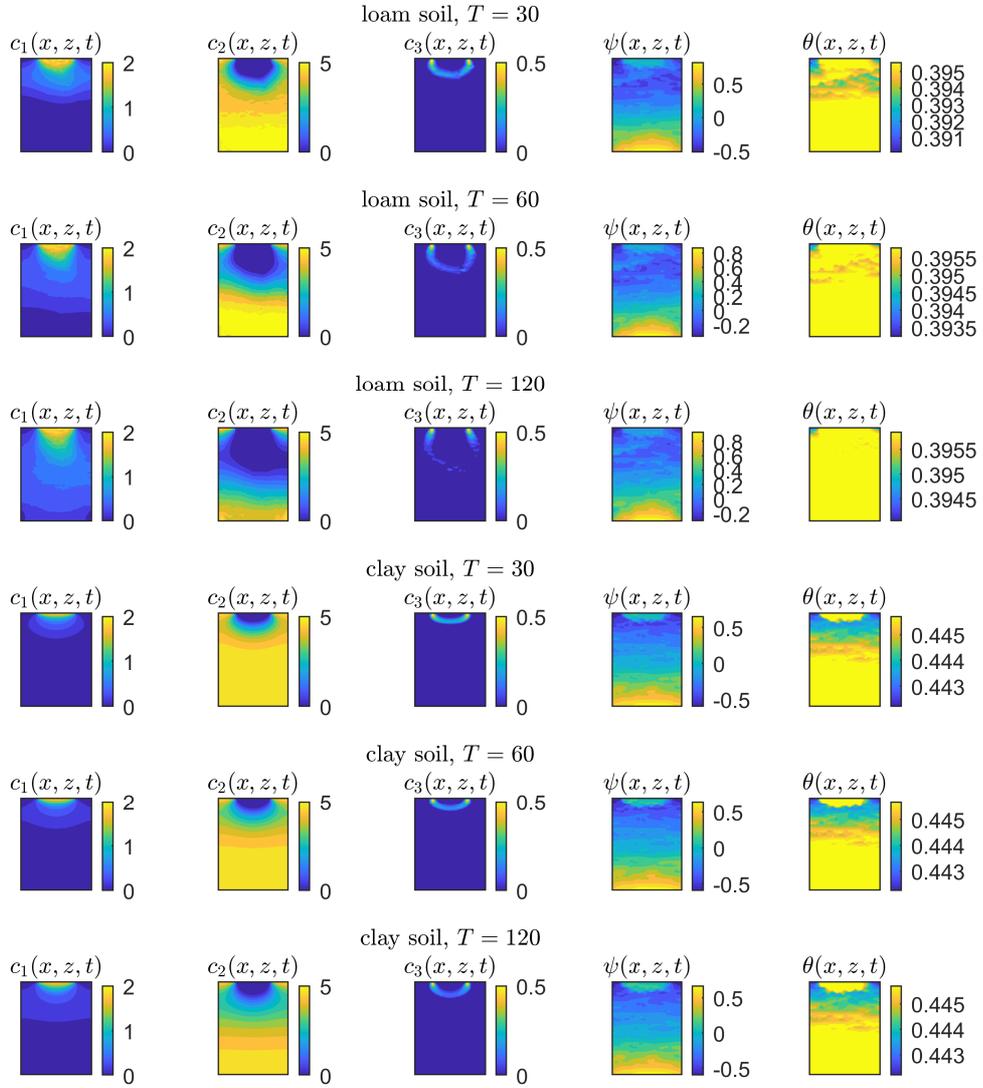}
\caption{\label{fig_loam_clay}Concentrations of the electron donor ($c_1$), electron acceptor ($c_2$), and biomass ($c_3$), pressure head ($\psi$), and water content ($\theta$) computed by the BGRW $L$-scheme coupled with the GRW flow solver for parameters of the loam and clay soil models.}
\end{figure}

The progress of the biodegradation and saturation processes for loam and clay soils is shown in Fig.~\ref{fig_loam_clay}. In the loam soil simulations the biodegradation process is quite fast, the oxygen is largely consumed, and the saturation is almost complete at $T=120$. In case of clay soil simulations, even though the pressure profiles are similar, only the lower part of the domain $\Omega$ becomes saturated and the biodegradation process is much slower. In addition, we observe that in case of loam soil the heterogeneity of the soil hydraulic properties influences not only the flow pattern but also the biodegradation process, as shown by the irregular contours of the benzene, oxygen, and biomass concentration. In clay soil simulations instead, even after 120 days the heterogeneity does not disturb significantly the regular concentration profiles, which are rather similar to those for homogeneous soils. These differences, consistent to those observed in similar simulations for nonreactive transport presented in \cite[Sect. 5.2.3]{Suciuetal2021}, can be explained not only by the lower hydraulic conductivity of the loam soil, but also by differences in $\theta_{sat}$, $\theta_{res}$, and parameters of the van Genuchten-Mualem model.

Based on the present simulations, different biodegradation scenarios can be constructed by slightly modifying the numerical setup, e.g., the injection conditions, the heterogeneity of the medium, as described by the random $K_{sat}$, the dispersion coefficients, or the Monod parameters.

\section{Conclusions}
\label{sec:concl}

The GRW algorithms for reactive transport belong to a larger class of particle methods for solving partial differential equations. These methods use random walkers or trajectories of the It\^{o} process to represent the solution and are, in principle, free of numerical diffusion. However, the GRW solvers differ in two important respects from other particle approaches. The first one is the global procedure which updates the positions of all the particles from a lattice site with computational costs that are practically the same as for a single random walk or PT step. The procedure allows us to use the same numbers of particles representing species concentrations as the number of molecules involved in chemical reactions. Consequently, there is no need to use spatial or ensemble averaging to obtain smooth concentration fields, the reaction steps are performed deterministically, and one avoids the cumbersome issue of representing the background oxygen concentration by particles filling the entire domain in PT simulations of biodegradation processes. The second advantageous feature is the robust coupling between the flow and reactive transport solvers, with flow solutions, provided by a GRW solver as well, computed at the sites of the lattice used in GRW transport simulations. One avoids thus the interpolation errors which occur when, for instance, the velocity field is computed with an external flow solver. This integrated approach is particularly useful in solving systems of coupled nonlinear flow and transport equations governing processes in saturated/unsaturated soils which require linearization procedures.

GRW schemes using unbiased random walk probabilities move large groups of particles from lattice sites by advective displacements and diffusive jumps, similarly to a PT scheme, with the difference that the particles undergo discrete movements between lattice sites, performed with a global procedure. BGRW schemes use random walk probabilities biased in the direction of the local velocity which simulate both the advection and the diffusion displacements through collective jumps on first-neighbor lattice sites. The two schemes can be used alternatively in solving specific transport problems. The unbiased GRW, which is significantly faster than the BGRW algorithm, is more efficient in solving large scale problems in unbounded domains but less efficient in solving boundary value problems. The BGRW instead is more accurate, allows straightforward implementations of the boundary conditions, and is more appropriate to solve coupled problems for soil systems with complex geometry.

Since the numbers of particles at lattice sites are readily converted into concentrations, we propose an operator splitting solution for the reaction step consisting of directly approximating the deterministic reaction terms and adding their contributions to the concentrations resulted from the advection-dispersion transport step. For saturated flow conditions the reactive transport solver can be implemented as a non-iterative scheme, with an explicit linearization of the reaction terms using concentration values provided by the transport step. Code verification tests for reactive transport with Monod reactions indicate the convergence towards manufactured analytical solutions of order one of the non-iterative GRW scheme and of order two of the BGRW scheme. For unsaturated flow conditions, the non-iterative schemes are no longer accurate and fail to converge (see Remark~\ref{rem:unsat-noniter}). Therefore, we designed iterative BGRW $L$-schemes coupled with GRW $L$-schemes for Richards equation, which can be easily adapted to bidirectional coupling between the flow and transport equations. The coupled $L$-schemes provide accurate solutions in case of variable saturation, as well as in case of degenerate Richards equation which describes transitions between unsaturated and saturated regimes. Code verification tests for transport and biodegradation indicate the approach to the convergence order 2 for pressure head and species concentrations solutions and a slower convergence of order 1.5 for the velocity components.

We applied the non-iterative GRW and BGRW schemes, which are faster than the iterative BGRW $L$-scheme, to solve a benchmark problem for biodegradation controlled by the transverse mixing in heterogeneous aquifers. The BGRW solution, for parameters describing a situation with moderate heterogeneity, served as reference to evaluate the faster but less accurate GRW solution. The comparison shows that the GRW solution reproduces the shape and the spatial distribution of electron donor, electron acceptor, and biomass concentrations but is affected by overshooting errors leading to the overestimation or the biodegradation process. Further investigations with the GRW solver show that the overshooting is enhanced by increasing the correlation lengths and the variance of the log-hydraulic conductivity. Since unbiased GRW and PT procedures move particles according to numerical schemes for It\^{o} equation over spatial regions with varying properties (e.g., reaction fronts), the overshooting can be only hardly avoided at the cost of using very small time steps (see Remark~\ref{rem:overshoot}). The occurrence of the overshooting errors in solving biodegradation problems remains a serious issue of concern for both unbiased GRW and PT approaches.

Finally, we solved the more challenging problem consisting of flow in soils, with transition from unsaturated to saturated conditions, and reactive transport with bidegradation governed by a double Monod model by using the flow GRW $L$-scheme coupled with the accurate BGRW $L$-scheme. The numerical solutions illustrate the progress of the saturation and biodegradation process in heterogeneous soils with different hydraulic properties and demonstrate the capability of the approach to simulate bioremediation scenarios.

\section*{Acknowledgements}
The first author acknowledges the financial support of the Deutsche Forschungsgemeinschaft (DFG, German Research Foundation) under Grant SU 415/4-1 -- 405338726 ``Integrated global random walk model for reactive transport in groundwater adapted to measurement spatio-temporal scales''.



\end{document}